\documentclass[11pt,a4paper,leqno]{article}
\usepackage[latin1]{inputenc}
\usepackage{amsmath}
\usepackage{amsfonts}
\usepackage{amssymb}

\newtheorem{satz}{Satz}[section]
\newtheorem{theorem}[satz]{Theorem}
\newtheorem{lemma}[satz]{Lemma}

\newtheorem{defin}[satz]{Definition}
\newtheorem{examples}[satz]{Examples}

\newcommand{\abs}[1]{\left|{#1}\right|}
\newcommand{\rund}[1]{\left(#1\right)}
\newcommand{\spitz}[1]{\left\langle{#1}\right\rangle}
\newcommand{\eckig}[1]{\left[{#1}\right]}
\newcommand{\schweif}[1]{\left\{#1\right\}}
\newcommand{\floor}[1]{\left\lfloor{#1}\right\rfloor}

\def\zz{\mathbb{Z}}
\def\cz{\mathbb{C}}
\def\nz{{\rm I\kern-.20em N}}
\def\rz{{\rm I\kern-.20em R}}
\def\pz{{\rm I\kern-.20em P}}

\def\C{\mathcal{C}}
\def\D{\mathcal{D}}
\def\E{\mathcal{E}}
\def\M{\mathcal{M}}
\def\N{\mathcal{N}}
\def\O{\mathcal{O}}
\def\F{\mathcal{F}}
\def\S{\mathcal{S}}
\def\T{\mathcal{T}}
\def\P{\mathcal{P}}
\def\Q{\mathcal{Q}}
\def\G{\mathcal{G}}
\def\cZ{\mathcal{Z}}

\def\a{\mathfrak a}
\def\g{\mathfrak g}
\def\z{\mathfrak z}
\def\u{\mathfrak u}
\def\sl{\mathfrak sl}
\def\su{\mathfrak su}
\def\m{\mathfrak m}
\def\n{\mathfrak n}

\def\vol{{\rm vol\ }}
\def\Im{{\rm Im\ }}
\def\hol{{\rm hol}}
\def\twist{{\rm twist}}
\def\Aut{{\rm Aut}}
\def\id{{\rm id}}
\def\Id{{\rm Id}}
\def\ad{{\rm ad}}
\def\Ad{{\rm Ad}}
\def\diag{{\rm diag}}
\def\nilp{{\rm nilp}}
\def\ord{{\rm ord}}
\def\tr{{\rm tr}}
\def\Ber{{\rm Ber \ }}
\def\Eig{{\rm Eig}}
\def\sdim{{\rm sdim\ }}

\def\bu{{\bf u} }
\def\bv{{\bf v} }

\def\Z{\begin{array}{c} z \\ \hline
\zeta
\end{array}}

\def\eps{\varepsilon}

\begin{document}

\vskip 1.0 true cm

\begin{center}
{\Huge \bf Super automorphic forms on}

\vskip 0.4 true cm

{\Huge \bf the super upper half plane}

\vskip 1.0 true cm

Roland Knevel\\

Unité de Recherche en Mathématiques Luxembourg

Campus Kirchberg
\end{center}

\vskip 2.0 true cm

\section*{Mathematical Subject Classification}

32C11 (Primary) , 11F55 (Secondary) .

\section*{Keywords}

Super symmetry, automorphic and cusp forms, local deformation of lattices, holomorphic vector bundles on compact Riemann surfaces.

\section*{Abstract}

Let $H^{|r}$ denote the upper half plane $H$ with $r$ additional odd (anticommuting) coordinates. It admits a transitive super action of a certain super Lie group $\G$ . First we define the spaces of super automorphic and cusp forms on $H^{|r}$ for an 
ordinary lattice $\Gamma$ of $\G$ , give an asymptotic formula for their dimensions for high weight and show how to embed $\Gamma \backslash H^{|r}$ into the super projective space with the help of super automorphic forms. For involving also the odd 
directions of $\G$ we introduce local super deformation of lattices in $\G$ and show that for high weight the spaces of super automorphic and cusp forms are stable under such local super deformations.

\section*{Introduction}

By now, super symmetry has been a current topic in physics for a long time with fruitfull influence on mathematics: Algebraic super structures and super manifolds were first invented as suitable mathematical tools for describing super symmetry in 
physics, but then they became more and more an independent field of research because of the elegance of the theory itself and the natural appearence among well-known classical mathematical structures, think for example of sheaves of differential forms. 
In purely mathematical context 'super' means: add 'odd (anticommuting) directions' to 'classical' objects. This leads to $\zz_2$-graded structures and the notion of super commutativity. So the theory of super manifolds embeds into 
the wide field of non-commutative geometry. I do not want to give a complete introduction to super manifolds here, the reader is referred to the literature, for example \cite{Ber}, \cite{Const} or \cite{KneBuch}. However, in section 1 I will briefly recall the definition and 
basic properties of the super upper half plane $H^{|r}$ as a super domain. Let me remark that there are two almost equivalent ways of describing super manifolds: via super numbers and via ringed spaces. Here I prefer the second one since it is more adapted to function spaces. \\

Riemann surfaces have continuously been objects of interest in mathematics. Most of them can be written as the upper half plane $H$ divided by a lattice in $\Aut H \simeq SL(2, \rz) / \{\pm 1\}$ , which leads to a relatively simple description of their moduli spaces, see \cite{Nat}. Finally automorphic forms play an important role in mathematics because of their connections to number theory, representation theory and algebraic geometry. For physicists they are of interest as an example of geometric quantization. In this article these three concepts will be combined. \\

We let a certain real super Lie group $\G$ act on $H^{|r}$ , and we want to fix a 'lattice' in $\G$ . A simple calculation shows that any $(0,0)$-dimensional sub super Lie group of $\G$ is nothing but an ordinary discrete subgroup in the body $G$ of $\G$ , and 
so up to this level we can forget about the odd directions of $\G$ . So how can we generalize the notion of a 'lattice' in $\G$ in order to involve also the odd directions of $\G$ ? The answer is: local super deformation. A single lattice in $\G$ has no chance to 
see the odd directions, but a whole family of lattices of course does if at least some of the 'parameters' parametrizing the family are odd and so all 'parameters' together generate a super commutative super algebra $\P$ . Such families will be called 
$\P$-lattices, they are local super deformations of the embedding of a single lattice into $\G$ . \\

There is some hope that as in the classical case super automorphic forms for a $\P$-lattice $\Upsilon$ will become a tool for 

\begin{itemize}
\item decomposing the left translation of the super Lie group $\G$ on some space of super functions on $\G / \Upsilon$ . The first aim will be to find an appropriate analogon for the classical $L^2$-space since integrability conditions do not make sense in 
the case of a $\P$-lattice,
\item identifying the quotient $\Upsilon \left\backslash H^{|r} \right.$ with some super algebraic variety, see theorem \ref{main} (iii) as a first step.
\end{itemize}

The paper is organized as follows: In section 2 we deal with the case of an 'ordinary' lattice in $G$ . Already this case is not at all trivial, and we give an asymptotic formula for the dimension of the spaces of super automorphic and super cusp forms for 
high weight $k$ , see theorem \ref{main} . This is done by writing super automorphic forms as global sections of vector bundles on the compact Riemann surface $X := \left.\Gamma^\# \right\backslash H \cup \{\text{cusps}\}$ , where $\Gamma^\#$ denotes the underlying lattice in $\Aut H$ . \\

While the classical deformation theory of lattices is already well-established, see \cite{Ragh}, in section 3 we talk about the generalization to the super case, giving both precise definitions, non-trivial examples and the connection with cohomology. 
Finally in section 5 we discuss super automorphic and cusp forms for $\P$-lattices. The main result here is the stability of the space of super automorphic forms for an ordinary lattice under its local super deformations for high weight $k$ , see theorem 
\ref{parammain}, which is obtained as a special case of local sheaf deformation discussed in section 4. In the special case $r = 0$ (so the usual upper half plane without odd coordinates) one already knows stability as soon as $k \geq 2$ or the genus $g$ of $\Gamma \backslash H \cup \{\text{cusps}\}$ is $\leq 1$ by a different method, see \cite{Kne} section 6. There one also finds a counterexample for the remaining case $k = 1$ and $g \geq 2$ . \\

{\it Acknowledgement:} I have to thank M. Schlichenmaier from Luxembourg and T. Bauer from Marburg for many helpful comments during the writing process and the Fonds National de la Recherche Luxembourg for funding my research stay at Luxembourg university.

\section{the general setting}

Let $r \in \nz$ (later in section 5 we have to exclude the case $r = 2$ ) and $\rund{GL(2, \cz) \times GL(r, \cz)}^{| 4 r}$ be the complex super Lie group with body $GL(2, \cz) \times GL(r, \cz)$ and $4 r$ additional odd (anticommuting) complex coordinates, 
where we sum up the $4 + r^2$ even and $4 r$ odd complex coordinates into an even super matrix

\[
g = \rund{\begin{array}{cc|c}
a & b & \mu \\
c & d & \nu \\ \hline
\rho & \sigma & E
\end{array}} \begin{array}{c} \left. \begin{array}{c}
 \\
\frac{}{}
\end{array} \right\} 2 \\
\left. \begin{array}{c}
 \\

\end{array} \right\} r
\end{array} \, .
\]

The equations $g I g^* = I$ and $\Ber g = 1$ , where $I := \rund{\begin{array}{c|c}
\begin{array}{cc}
0 & i \\
- i & 0
\end{array} & 0 \\ \hline
0 & 1
\end{array}}$ , and $\Ber g := \det\rund{\rund{\begin{array}{cc}
a & b \\
c & d
\end{array}} - \rund{\begin{array}{c}
\mu \\
\nu
\end{array}} \rund{\begin{array}{cc}
\rho & \sigma
\end{array}} } \det E^{- 1}$ denotes the super determinant (the so-called Berezinian) of $g$ , define a real super Lie group $\G$ of super dimension $(3 + r^2, 4 r)$ with body

\[
G := \schweif{\left.\rund{\begin{array}{c|c}
\eps h & 0 \\ \hline
0 & E
\end{array}} \, \right| \, \eps \in U(1) , h \in SL(2, \rz) , E \in U(r) , \eps^2 = \det E}
\]

and super Lie algebra $\g = \g_0 \oplus \g_1$ ,

\begin{eqnarray*}
\g_0 &:=& \schweif{\left.\rund{\begin{array}{c|c}
\begin{array}{cc}
a + \frac{1}{2} \tr D & b \\
c & - a + \frac{1}{2} \tr D
\end{array} & 0 \\ \hline
0 & D
\end{array}} \, \right| \, a, b, c \in \rz , D \in \u(r)} \\
&& \phantom{\frac{1}{23}} \simeq \sl(2, \rz) \oplus \u(r) \, , \\
\g_1 &:=& \schweif{\left.\rund{\begin{array}{c|c}
0 & \begin{array}{c}
\bv^* \\
- \bu^*
\end{array} \\ \hline
\begin{array}{cc}
\bu & \bv \end{array} & 0
\end{array}} \, \right| \, \bu, \bv \in \cz^r} \, .
\end{eqnarray*}

Let $H^{|r}$ denote the usual upper half plane $H := \{\Im > 0\} \subset \cz$ with $r$ additional odd (anticommuting) complex coordinate functions. Then $H^{|r}$ is in particular a super domain, and we recall the basic properties. \\

As a super domain $H^{|r}$ is defined as the ringed space $\rund{H, \O_H \otimes \bigwedge\rund{\cz^r}}$~, $\O_H \otimes~\bigwedge\rund{\cz^r}$ being a sheaf of complex unital associative super commutative super algebras, the sheaf of holomorphic 
super functions on $H^{|r}$ (by definition $\zz_2$-graded). This sheaf even admits a $\zz$-grading coming from the well-known $\zz$-grading $\bigwedge\rund{\cz^r} = \bigoplus_{\rho = 0}^r \bigwedge^\rho\rund{\cz^r}$ of the exterior algebra, and for 
$U \subset H$ open we write $\O\rund{U^{|r}} := \O(U) \otimes \bigwedge\rund{\cz^r} = \bigoplus_{\rho = 0}^r \O^\rho\rund{U^{|r}}$ , where $\O^\rho\rund{U^{|r}} := \O(U) \otimes \bigwedge^\rho\rund{\cz^r}$ . The odd complex coordinates of $H^{|r}$~, which 
are nothing but the standard basis vectors in $\bigwedge^1\rund{\cz^r} = \cz^r$ , will always be denoted by $\zeta_1, \dots, \zeta_r$ . We denote the power set of $\{1, \dots, r\}$ by $\wp(r)$ , and for every $I \in \wp(r)$ , $I = \schweif{i_1, \dots, i_\rho}$ , 
$i_1 < \dots < i_\rho$ , we write

\[
\zeta^I := \zeta_{i_1} \cdots \zeta_{i_\rho} \, .
\]

Therefore every holomorphic super function $f \in \O\rund{U^{|r}}$ has a unique decomposition $f = \sum_{I \in \wp(r)} f_I \zeta^I$ , all $f_I \in \O(U)$ . The super automorphisms of $H^{|r}$ are by definition the
automorphisms of $H^{|r}$ as a ringed space. So every super automorphism $\Phi$ of $H^{|r}$ has an underlying ordinary automorphism $\Phi^\# \in \Aut H$ , which is called the body of $\Phi$ . In practice the super automorphisms of 
$H^{|r}$ are given by tuples $\rund{f, \lambda_1, \dots, \lambda_r} \in \O\rund{H^{|r}}_0 \oplus \O\rund{H^{|r}}_1^{\oplus r}$~, and in this notation the body is given by $f^\# \in \O(H)$ , where we denote by 
${}^\#$ the $\O_H$-linear extension of the canonical projection ${}^\#: \bigwedge\rund{\cz^r} \rightarrow \cz$ . \\

We have a transitive holomorphic super action $\alpha: \G \times H^{|r} \rightarrow H^{|r}$ of $\G$ on $H^{|r}$ given by super Möbius transformations

\[
g \rund{\Z} := \frac{1}{c z + d + \nu \zeta} \rund{\begin{array}{c}
a z + b + \mu \zeta \\ \hline
\rho z + \sigma + E \zeta
\end{array}} \, , \, g = \rund{\begin{array}{cc|c}
a & b & \mu \\
c & d & \nu \\ \hline
\rho & \sigma & E
\end{array}} \, .
\]

Its body $\alpha^\#: G \otimes H \rightarrow H$ extends the well-known action of 

\[
SL(2, \rz) \hookrightarrow G \, , \, h \mapsto \rund{\begin{array}{c|c}
h & 0 \\ \hline
0 & 1
\end{array}}
\]

on $H$ by classical Möbius transformations. By $\alpha$ we have a group homomorphism from $G$ into the group of super automorphisms of $H^{|r}$ , and if we apply in addition the body functor from $H^{|r}$ to $H$ to these super automorphisms we 
even obtain a group homomorphism

\[
{}^\# : G \rightarrow \Aut H \simeq SL(2, \rz) / \{\pm 1\} \, , \, \rund{\begin{array}{c|c}
\eps h & 0 \\ \hline
0 & E
\end{array}} \mapsto \overline h \, ,
\]

$\eps \in U(1)$ , $h \in SL(2, \rz)$ , $E \in U(r)$ , $\eps^2 = \det E$ . $G_0 := \ker {}^\# \sqsubset G$ is a compact subgroup. Since $G$ is an almost direct product of $SL(2, \rz)$ and $G_0 \sqsubset G$ we see that $G$ is unimodular. \\

By the way, via a super Cayley transform mapping biholomorphically the super unit disc $B^{|r}$ onto the super upper half plane $H^{|r}$ by super Möbius transform this situation is equivalent to the one treated in \cite{Borth} , where the super Lie group 
$SU(1, 1|r)$ acts on $B^{|r}$ via super Möbius transformations. \\

For a lattice $\Gamma \sqsubset G$ , which means by definition discrete of finite covolume, we define $\Gamma_0 := \Gamma \cap G_0 \sqsubset G_0$ finite , $\Gamma^\# := \schweif{\left.\gamma^\# \, \right| \, \gamma \in \Gamma} \sqsubset \Aut H$ and 
$\check \Gamma$ to be the preimage of $\Gamma^\#$ under the canonical projection $\overline{\phantom{1}} : SL(2, \rz) \rightarrow \Aut H$~. Then $\check \Gamma \sqsubset SL(2, \rz) \hookrightarrow G$ is at the same time the set of all $h \in SL(2, \rz)$ 
such that there exists $\eta \in G_0$ with $h \eta \in \Gamma$ . Moreover:

\begin{lemma}
$\Gamma^\# \sqsubset \Aut H$ and $\check \Gamma \sqsubset SL(2, \rz)$ are lattices.
\end{lemma}

{\it Proof:} $\Gamma^\#$ and $\check \Gamma$ are trivially discrete. For proving that $\Gamma^\#$ and $\check \Gamma$ are of finite covolume let $\Omega \subset SL(2, \rz)$ be open such that $\check \Gamma \, \Omega = SL(2, \rz)$ and 
$\check \gamma \Omega \cap \Omega \not= \emptyset$ for only finitely many $\check \gamma \in \check \Gamma$ . Then the same is true for $\Omega \, G_0$ with respect to $\Gamma$ . So $\vol \Omega \, G_0 < \infty$ and so also $\vol \Omega < \infty$ . 
$\Box$ \\

So $X := \left.\Gamma^\# \right\backslash H \cup \schweif{\text{cusps of } \left.\Gamma^\# \right\backslash H}$ has the structure of a compact Riemann surface. Let $\pi_X : H \rightarrow \left.\Gamma^\# \right\backslash H \hookrightarrow X$ denote the 
canonical projection. Let $z_0 \in \partial_{\pz^1} H$ , then there exists $g \in \Aut H$ such that $g \, i \infty = z_0$~. For using the standard notation we call $N^{z_0} := g N^{i \infty} g^{- 1} \sqsubset \Aut H$ the nilpotent subgroup associated to $z_0$ , where 
$N^{i \infty}$ is the image of the group embedding $\rz \hookrightarrow \Aut H$ assigning to $t \in \rz$ the translation $z \mapsto z + t$ , and we call an open set $U \subset H$ a neighbourhood of $z_0$ iff there exists $R > 0$ such that 
$g \, \{\Im z > R\} \subset U$ . If $z_0$ is a cusp of $\left.\Gamma^\# \right\backslash H$ then the neighbourhoods of $z_0$ in $H$ are precisely the subsets $U \subset H$ such that $\pi_X(U)$ is a punctured neighbourhood of $\overline{z_0}$ in $X$ . \\

In the end of this section let us discuss two examples of lattices $\Gamma \sqsubset G$ : \\

\begin{examples} \label{exordlat} \end{examples}

Let $\gamma_0 := \rund{\begin{array}{c|c}
\eps_0 1_2 & 0 \\ \hline
0 & E_0
\end{array}} \in G_0$ be of finite order $N$ with $E_0 \in U(r)$ , \\
$\eps_0 \in U(1)$ , $\eps_0^2 = \det E_0$ .

\begin{itemize}
\item[$\spitz{\text{i}}$] $R := \rund{\begin{array}{cc}
0 & 1 \\
- 1 & - 1
\end{array}}$ and $S := \rund{\begin{array}{cc}
0 & 1 \\
- 1 & 0
\end{array}} \in SL(2, \rz)$ generate $SL(2, \zz)$ . Let furthermore $E, F \in Z_{U(r)}\rund{E_0}$ such that $E^3 = E_0^m$~, $F^2 = E_0^n$ and $\det F = - \eps_0^n$ for some $m, n \in \nz$ , $\eps, \eta \in U(1)$ such that $\eps^2 = \det E$ , $\eta^2 = \det F$ and 
$\eps^3 = \eps_0^m$ . Then $\gamma_0$ , 

\[
\hat R := \rund{\begin{array}{c|c}
\eps R & 0 \\ \hline
0 & E
\end{array}} \phantom{12} \text{ and } \phantom{12} \hat S := \rund{\begin{array}{c|c}
\eta S & 0 \\ \hline
0 & F
\end{array}}
\]

generate a lattice $\Gamma \sqsubset G$ with $\Gamma_0 = \spitz{\gamma_0}$ and $\check \Gamma = SL(2, \zz)$ . It is the free group in $\gamma_0$ , $\hat R$ and $\hat S$ moludo the relations

\[
{\hat R}^3 = \gamma_0^m \, , \, \hat S^2 = \gamma_0^n \, , \, \eckig{\hat R, \gamma_0} = \eckig{\hat S, \gamma_0} = \gamma_0^N = 1 \, .
\]

\item[$\spitz{\text{ii}}$] Let $X$ be a compact Riemann surface of genus $g^*$ , $m \in \nz$ , $3 g^* + m \geq~4$~, and $s_1, \dots, s_m \in X$ . Then the universal covering of $X \setminus \schweif{s_1, \dots, s_m}$ is isomorphic to $H$ , and by \cite{Nat} 
one can write $X \setminus \schweif{s_1, \dots, s_m} = \left.\Gamma' \right\backslash H$~, where $\Gamma' \sqsubset SL(2, \rz)$ is a lattice without elliptic elements, $- 1 \notin \Gamma'$ and 
$\Gamma' \simeq \pi_1 \rund{X \setminus \schweif{s_1, \dots, s_m}}$ . It is the free group generated by some hyperbolic elements $A_1, B_1, \dots, A_{g^*}, B_{g^*}$ and parabolic elements $C_1, \dots, C_m \in SL(2, \rz)$ modulo the single relation 
$\eckig{A_1, B_1} \cdots \eckig{A_{g^*}, B_{g^*}} C_1 \cdots C_m = 1$ .

Let furthermore $E_k, F_k, H_l \in Z_{U(r)}$ , $k = 1, \dots, g^*$ , $l = 1, \dots m$ , such that $\eckig{E_1, F_1} \cdots \eckig{E_{g^*}, F_{g^*}} H_1 \cdots H_m = E_0$ and $\eps_k, \eta_k, \vartheta_l \in U(1)$ such that $\eps_k^2 = \det E_k$ , 
$\eta_k^2 = \det F_k$ , $\vartheta_l^2 = H_l$ , $k = 1, \dots, g^*$ , $l = 1, \dots m$ , and $\vartheta_1 \cdots \vartheta_m = \eps_0$ . Then $\gamma_0$ ,

\begin{eqnarray*}
&& \hat A_k := \rund{\begin{array}{c|c}
\eps_k A_k & 0 \\ \hline
0 & E_k
\end{array}} \, , \, \hat B_k := \rund{\begin{array}{c|c}
\eta_k B_k & 0 \\ \hline
0 & F_k
\end{array}} \, , \, k = 1, \dots, g^* \, , \\
&& \phantom{12} \text{ and } \phantom{12} \hat C_l := \rund{\begin{array}{c|c}
\vartheta_l C_l & 0 \\ \hline
0 & H_l
\end{array}} \, , l = 1, \dots, m \, ,
\end{eqnarray*}

generate a lattice $\Gamma \sqsubset G$ with $\Gamma_0 = \spitz{\gamma_0}$ , $\Gamma^\# = \overline{\Gamma'}$ and $\eckig{\hat A_1, \hat B_1} \cdots \eckig{\hat A_{g^*}, \hat B_{g^*}} \hat C_1 \cdots \hat C_m = \gamma_0$ . \\

If $m \geq 1$ then $\Gamma$ is the free group in the generators $\gamma_0$ , $\hat A_k$ , $\hat B_k$ , $k = 1, \dots, g^*$ , and $\hat C_l$ , $l = 1, \dots, m - 1$ (!) , modulo the relations

\[
\eckig{\hat A_k, \gamma_0} = \eckig{\hat B_k, \gamma_0} = \eckig{\hat C_l, \gamma_0} = \gamma_0^N = 1 \, .
\]

If $m = 0$ then necessarily $\gamma_0 = 1$ , and so $\Gamma$ is the free group in the generators $\hat A_k \, , \, \hat B_k$ , $k = 1, \dots, g^*$ , moludo the single relation

\[
\eckig{\hat A_1, \hat B_1} \cdots \eckig{\hat A_{g^*}, \hat B_{g^*}} = 1 \, .
\]

\end{itemize}

\section{super automorphic forms for ordinary lattices}

On $\G \times H^{|r}$ we have a cocycle $j \in \rund{\D(\G)^\cz \hat\boxtimes \O\rund{H^{|r}} }_0$ , where \\
$\D(\G) \simeq \C^\infty(G) \otimes \bigwedge\rund{\rz^{4 r}}$ , $G$ being the body of $\G$ , denotes the space of (real valued smooth) super functions on $\G$ and by '$\boxtimes$' we denote the $\zz_2$-graded tensor product, given by

\[
j\rund{g, \Z} := \frac{1}{c z + d + \nu \zeta} \phantom{12} , \phantom{12} g = \rund{\begin{array}{cc|c}
a & b & \mu \\
c & d & \nu \\ \hline
\rho & \sigma & E
\end{array}} \, ,
\]

and for each $k \in \zz$ the assignment 

\[
f \mapsto f\rund{g \rund{\Z}} \, j\rund{g, \Z}^k
\]

defines a $\zz_2$-graded linear map $|_k : \O\rund{H^{|r}} \rightarrow \D(\G)^\cz \hat\boxtimes \O\rund{H^{|r}}$ and for each $g \in G$ and $U \subset H$ open a $\zz$-graded (!) linear map $|_{g, k} : \O\rund{\rund{g^\# U}^{|r}} \rightarrow \O\rund{U^{|r}}$ . 
Usually we will drop the index $k$ . Observe that the Berezinian of the super Jacobian of $\alpha$ with respect to $\rund{\Z}$ is precisely given by $j^{2 - r}$ , see \cite{Borth}. \\

For defining super automorphic resp. cusp forms for a lattice $\Gamma$ first we have to give a notion of boundedness resp. vanishing of a super function on the super upper half plane $H^{|r}$ at a cusp of $\left.\Gamma^\# \right\backslash H$ . For this 
purpose let 

\begin{equation} \label{standardg_0}
g_0 := \rund{\begin{array}{c|c}
\eps_0 \rund{\begin{array}{cc}
1 & 1 \\
0 & 1
\end{array}} & 0 \\ \hline
0 & E_0
\end{array}} \in G \, ,
\end{equation}

$\eps_0 \in U(1)$ , $E_0 = \rund{\begin{array}{ccc}
e_1 & & 0 \\
 & \ddots & \\
0 & & e_r
\end{array}}\in U(r)$ diagonal, $\eps_0^2 = \det E_0$ , and let $f = \sum_{I \in \wp(r)} f_I \zeta^I \in \O\rund{\{\Im z > R\}^{|r}}$ , $R > 0$ . Then for all $I \in \wp(r)$

\[
\left.f_I \zeta^I \right|_{g_0} = f_I(z + 1) \, \eps_0^{- k - \abs{I}} \det\nolimits_I E_0 \, \zeta^I \, ,
\]

where $\det_I E_0 := \prod_{i \in I} e_i$ , and $\eps_0^{- k - \abs{I}} \det_I E_0 \in~U(1)$~. So if $f|_{g_0} = f$ then all $f_I$ are quasi-invariant under $z \mapsto z + 1$ .

\begin{defin} \label{bound} \end{defin}

\begin{itemize}

\item[(i)] Let $R > 0$ and $f = \sum_{I \in \wp(r)} f_I \zeta^I \in \O\rund{\{\Im z > R\}^{|r}}$ such that \\
$f|_{g_0} = f$~. Then $f$ is called bounded (vanishing) at $i \infty$ iff all $f_I(z)$ , $I \in \wp(r)$ , are bounded (vanishing) for $\Im z \leadsto \infty$ .
\item[(ii)] Let $z_0 \in \partial_{\pz^1} H$ and $\gamma \in G$ such that $\gamma^\# \in N^{z_0} \setminus \{\id\}$ . Let $U \subset H$ be an open $\gamma^\#$-invariant neighbourhood of $z_0$ and $f \in \O\rund{U^{|r}}$ such that $f|_\gamma = f$ . Take 
some $g \in G$ such that $g^\# i \infty = z_0$ and either $g_0 := g^{- 1} \gamma g$ or $g_0 := g^{- 1} \gamma^{- 1} g$ is of the form (\ref{standardg_0}). Then $f|_g$ is invariant under $|_{g_0}$ . $f$ is called bounded (vanishing) at $z_0$ iff $f|_g$ is 
bounded (vanishing) at $i \infty$ .

\end{itemize}

Of course we have to prove invariance of definition \ref{bound} (ii) under the choice of $g \in G$ :

\begin{quote}
Let $g \in G$ such that $g_0' := g^{- 1} g_0 g$ is again of the form (\ref{standardg_0}) with some $\eps_0' \in U(1)$ , $E_0' \in U(r)$ diagonal, $\eps_0'^2 = \det E_0'$  . Then $\eps_0' = \eps_0$ , and so $g$ and $g_0$ commute in the upper left corner. 
Therefore

\[
g = \rund{\begin{array}{c|c}
\eps \rund{\begin{array}{cc}
1 & S \\
0 & 1
\end{array}} & 0 \\ \hline
0 & E
\end{array}}\]

with some $\eps \in U(1)$ , $S \in \rz$ , $E \in U(r)$ , $\eps^2 = \det E$ . So $f|_g$ is a linear combination of terms $f_I(z + S) \zeta^J$ , $I, J \in \wp(r)$ . $\Box$ \\
\end{quote}

Now let $\Gamma \sqsubset G$ be a lattice and $k \in \zz$ .

\begin{defin}[super automorphic and super cusp forms for $\Gamma$ ] \end{defin}

Let $f \in \O\rund{H^{|r}}$ . $f$ is called a super automorphic (cusp) form for $\Gamma$ of weight $k$ iff

\begin{itemize}
\item[(i)] $f|_{\gamma, k} = f$ for all $\gamma \in \Gamma$ ,
\item[(ii)] $f$ is bounded (vanishing) at all cusps of $\left.\Gamma^\# \right\backslash H$ in the sense of definition \ref{bound}.
\end{itemize}

The $\cz$- vector space of super automorphic (cusp) forms for $\Gamma$ of weight $k$ is denoted by $sM_k(\Gamma)$ (resp. $sS_k(\Gamma) \sqsubset sM_k(\Gamma)$ ). \\

Since $|_g$ respects the $\zz$-grading of $\O\rund{H^{|r}}$ for all $g \in G$ we obtain a $\zz$-grading $sM_k(\Gamma) = \bigoplus_{\rho = 0}^r sM_k^\rho(\Gamma)$ (resp. $sS_k(\Gamma) = \bigoplus_{\rho = 0}^r sS_k^\rho(\Gamma)$ ) where 
$sM_k^\rho(\Gamma) = sM_k(\Gamma) \cap \O^\rho\rund{H^{|r}}$ (resp. $sS_k^\rho(\Gamma) = sS_k(\Gamma) \cap \O^\rho\rund{H^{|r}}$ ). \\

\begin{examples} \label{exautfo} \end{examples}

\begin{itemize}

\item[$\spitz{\text{i}}$] If $\Gamma \sqsubset SL(2, \rz) \hookrightarrow G$ is a lattice then $sM_k^\rho(\Gamma) = M_{k + \rho}(\Gamma) \otimes \bigwedge^\rho\rund{\cz^r}$ for all $k \in \zz$ and 
$\rho = 0, \dots, r$ , where $M_{k + \rho}(\Gamma) \sqsubset \O(H)$ denotes the space of ordinary automorphic forms for $\Gamma$ of weight $k + \rho$ . In particular if $- 1 \in \Gamma$ and 
$k + \rho$ is odd then $sM_k^\rho(\Gamma) = 0$ . This behaviour corresponds to the philosophy of super symmetry to regard different sorts of bosons and fermions as components of one 
super particle. \\

\item[$\spitz{\text{ii}}$] As a special case of example \ref{exordlat} $\spitz{\text{i}}$ let $r = 1$ and $\Gamma \sqsubset G$ be the lattice generated by

\[
\hat R := \rund{\begin{array}{c|c}
\eps R & 0 \\ \hline
0 & \eps^{- 1}
\end{array}} \, \text{ and } \, \hat S := \rund{\begin{array}{c|c}
i \omega S & 0 \\ \hline
0 & - 1
\end{array}} \, .
\]

Let $\eta := e^{\frac{\pi i}{12} z} \prod_{n = 1}^\infty \rund{1 - e^{2 \pi i n z}}$ denote Dedekind's eta function. \\

If $\eps = e^{\frac{2 \pi i}{3}}$ and $\omega = 1$ then $\eta^2$ generates $sM_1^0(\Gamma)$ , and \\
$sM_1^1(\Gamma) = M_2(SL(2, \zz)) \zeta = 0$ .

If $\eps = e^{- \frac{2 \pi i}{3}}$ and $\omega = - 1$ then $\eta^2 \zeta$ generates $sM_0^1(\Gamma)$ , and $sM_0^0(\Gamma) = \cz$~. \\

In both cases we have $\Gamma_0 = \{1\}$ . The result can be seen by computing directly the vector bundles $E^\rho_k$ (in fact line bundles here) on $X$ (the Riemann sphere here), which 
will be defined in $\{\text{i}\}$, and their degrees: In the first case $\deg E_1^0 = 0$ and $\deg E_1^1 = - 1$ . In the second case $\deg E_0^1 = 0$ and $E^0_0$ is of course trivial. \\

Similarly also the square of the theta function $\theta := \sum_{n = - \infty}^\infty e^{\pi i n^2 z}$ can be realized as an even or odd super automorphic form.

\item[$\spitz{\text{iii}}$] Let $\Gamma \sqsubset G$ be a lattice such that $\overline \Gamma \sqsubset \overline{SL(2, \zz)}$ is of finite index and

\[
q := \min \schweif{q \in \nz \setminus \{0\}\, \left| \, \rund{z \mapsto z + q} \in \overline \Gamma\right.} \, .
\]

Then there exist $\eps_0 \in U(1)$ and $E_0 \in U(r)$ , which we may assume to be diagonal, $\eps_0^2 = \det E_0$ , such that

\[
\gamma_0 := \rund{\begin{array}{c|c} \eps_0 \rund{\begin{array}{cc}
1 & q \\
0 & 1
\end{array}} & 0 \\ \hline
0 & E_0
\end{array}} \in \Gamma \, .
\]

Let $\nu \in \rz_{\geq 0}$ such that $e^{2 \pi i q \nu} = \eps_0^{k + \abs{I}} \det_I E_0$ . Then the function $f = e^{2 \pi i \nu z} \zeta^I$ is already invariant under $|_{\gamma_0, k}$ , and a 
simple estimate shows that the relative super Poincaré series

\[
\sum_{\gamma \in \left.\spitz{\gamma_0} \right\backslash \Gamma} f|_{\gamma, k}
\]

converges absolutely and uniformly on compact sets of $H$ and defines an element in $sM_k(\Gamma)$ , lying in $sS_k(\Gamma)$ iff $\nu > 0$ .

\end{itemize}

Let $\rho \in \{0, \dots, r\}$ . Then for all $g \in G_0$ and $U \subset H$ open $|_g$ is an $\O(U)$-linear (!) operator on $\O^\rho\rund{U^{|r}}$ with $\bigwedge^\rho\rund{\cz^r}$ as invariant subspace. Moreover $|$ defines a unitary right 
representation of $G_0$ on $\bigwedge^\rho\rund{\cz^r}$~. Write

\[
V_k^\rho := \schweif{\left.a \in \bigwedge\nolimits^\rho\rund{\cz^r} \, \right| \, a|_{\eta, k} = a \text{ for all } \eta \in \Gamma_0} \sqsubset \bigwedge\nolimits^\rho\rund{\cz^r}
\]

for all $k \in \zz$ . Then obviously $sM_k^\rho(\Gamma) \sqsubset \O(H) \otimes V_k^\rho$ for all $k \in \zz$ . As a first observation we remark:

\begin{lemma}
The families $\rund{V_k^\rho}_{k \in \zz}$ of subspaces of $\bigwedge^\rho\rund{\cz^r}$ are $\abs{\Gamma_0}$-periodic.
\end{lemma}

{\it Proof:} Let $\eta \in \Gamma_0$ , which is then of the form $\eta = \rund{\begin{array}{c|c} \eps 1 & 0 \\ \hline
0 & E
\end{array}}$ with some $\eps \in U(1)$ and $E \in U(r)$ , $\eps^2 = \det E$ . Since $\eta^{\abs{\Gamma_0}} = 1$ we obtain $\eps^{\abs{\Gamma_0}} = 1$ . Now for all $\rho \in \{0, \dots, r\}$ , $a = \sum_{\abs{I} = \rho} a_I \zeta^I \in \bigwedge^\rho\rund{\cz^r}$ 
and $k \in \zz$ :

\[
a|_{\eta, k + \abs{\Gamma_0}} = \eps^{- \rund{\abs{\Gamma_0} + k + \rho}} \sum_{\abs{I} = \rho} a_I (E \zeta)^I = \eps^{- (k + \rho)} \sum_{\abs{I} = \rho} a_I (E \zeta)^I = a|_{\eta, k} \, . \, \Box
\]

Here now the main theorem, whose proof will be the purpose of the rest of this section:

\begin{theorem}[main theorem] \label{main} Let $\rho \in \{0, \dots, r\}$ .
\item[(i)] There exists $k_0 \in \zz$ such that $sM_k^\rho(\Gamma) = 0$ for all $k \leq k_0$ .

\item[(ii)] For $k \leadsto + \infty$ we have the asymptotic behaviour

\[
\dim sM_k^\rho(\Gamma) = \rund{\frac{k}{2} \, \frac{\vol\rund{\left.\Gamma^\# \right\backslash H}}{2 \pi} + O(1)} \dim V_k^\rho \, ,
\]

and for all $k \in \zz$

\[
\dim sS_k^\rho(\Gamma) \geq \dim sM_k^\rho(\Gamma) - S \dim V_k^\rho \, ,
\]

where $S$ denotes the number of cusps of $\left.\Gamma^\# \right\backslash H$ .

\item[(iii)] If $\Gamma_0 = \{1\}$ and $\left.{\Gamma^\#} \right\backslash H$ has no elliptic points then there exists \\
$k_1 \in \zz$ such that for all $k \geq k_1$ , given bases $\schweif{f_0, \dots, f_m}$ of $sM_k^0(\Upsilon)$ and 
$\schweif{\lambda_1, \dots, \lambda_n}$ of $sM_k^1(\Gamma)$ , $\Phi := \eckig{f_0: \dots: f_m: \lambda_1: \dots: \lambda_n}$ defines an embedding of $\Gamma \left\backslash H^{|r}\right.$ into the $(m, n)$-dimensional complex super projective 
space $\pz^{m|n}$ as a complex $(1, r)$-dimensional sub super manifold, if in addition $\left.{\Gamma^\#} \right\backslash H$ has no cusps then in fact as a non-singular algebraic sub super variety.
\end{theorem}

\begin{lemma} \label{repres}

\item[(i)] For all $k \in \zz$ and $\rho \in \{0, \dots, r\}$ there exists a unique unitary right representation $\varphi^\rho_k$ of $\check \Gamma$ on $V_k^\rho$ such that

\[
f|_\gamma = j\rund{\check \gamma, z}^{k + \rho} \, f\rund{\gamma^\# z} \varphi^\rho_k\rund{\check \gamma} \, ,
\]

for all $U \subset H$ open, $f \in \O(U) \otimes V_k^\rho \sqsubset \O^{\rho}\rund{U^{|r}}$ , $\gamma \in \Gamma$ and $\check \gamma \in \check \Gamma$ 
representing $\gamma^\#$ , where we have extended $\varphi^\rho_k\rund{\check \gamma}$ as an $\O(U)$-linear map to $\O(U) \otimes V_k^\rho$. 
Obviously $\varphi^\rho_k(- 1) = (- 1)^{k + \rho}$ and $\varphi_0^0 = 1$ .

\item[(ii)] Let $U \subset H$ be open and $\Gamma^\#$-invariant and $f \in \O^{\rho}\rund{U^{|r}}$. Then in particular $f$ is invariant under all $|_\gamma$ , $\gamma \in \Gamma$ , iff $f \in \O(U) \otimes V_k^\rho$ and

\[
j\rund{\check \gamma, z}^{k + \rho} \, f\rund{\check \gamma z} \varphi^\rho_k\rund{\check \gamma} = f
\]

for all $\check \gamma \in \check \Gamma$ .

\item[(iii)] There exists a unique character $\chi: \check \Gamma \rightarrow U(1)$ such that $\varphi^\rho_{k + \abs{\Gamma_0}} = \chi \cdot \varphi^\rho_k$ 
for all $k \in \zz$ and $\rho \in \{1, \dots, r\}$ .

\end{lemma}

{\it Proof:} (i) Let $\check \gamma \in \check \Gamma$ and $\eta \in G_0$ such that $\gamma := \eta \check \gamma \in \Gamma$ . Then

\[
f|_\gamma = \left.f|_\eta \right|_{\check \gamma} = j\rund{\check \gamma, z}^{k + \rho} \, f|_\eta \rund{\gamma^\# z} \, .
\]

So the formula gives the right expression for $f|_\gamma$ iff we define $a \ \varphi^\rho_k\rund{\check \gamma} := a|_\eta$ for all $a \in V_k^\rho$ . But we have to check that this indeed defines a unitary representation of $\check \Gamma$ 
on $V_k^\rho$ . \\

For proving well-definedness first we show that again $a|_\eta \in V_k^\rho$ : Let $\vartheta \in \Gamma_0$~. Then since $\Gamma_0 \lhd \Gamma$ again 
$\eta \vartheta \eta^{- 1} = \gamma \vartheta \gamma^{- 1} \in \Gamma_0$ . So

\[
\left.a|_\eta \ \right|_\vartheta = \left.a|_{\eta \vartheta \eta^{- 1}} \right|_\eta = a|_\eta \, .
\]

Now let us show that $a|_{\eta^{- 1}}$ is independent of the particular choice of $\eta$ : Let also $\vartheta \in G_0$ such that $\vartheta \check \gamma \in \Gamma$ . Then 
$\vartheta \eta^{- 1} = \rund{\vartheta \check \gamma} \gamma^{- 1} \in \Gamma_0$ . Therefore

\[
a|_\vartheta = \left.a|_{\vartheta \eta^{- 1}} \right|_\eta = a|_\eta \, .
\]

$\check \Gamma$ and $G_0$ commuting shows that $\varphi_k^\rho$ is indeed a right representation of $\check \Gamma$ .

(ii) now trivial.

(iii)  Let $\check \gamma \in \check \Gamma$ and $\eta = \rund{\begin{array}{c|c} \eps 1 & 0 \\ \hline
0 & E
\end{array}} \in~G_0$~, $\eps \in U(1)$ and $E \in U(r)$ , $\eps^2 = \det E$ , such that $\gamma := \eta \check \gamma \in \Gamma$ , and let $a = \sum_{\abs{I} = \rho} a_I \zeta^I \in V_k^\rho$ . Then by the 
proof of (i) we see that

\[
a \ \varphi^\rho_{k + \abs{\Gamma_0}}\rund{\check \gamma} = a|_{\eta, k + \abs{\Gamma_0}} = \eps^{- k - \abs{\Gamma_0} - \rho} \sum_{\abs{I} = \rho} a_I \rund{E \ \zeta}^I \, .
\]

So the formula holds iff we define $\chi\rund{\check \gamma} := \eps^{- \abs{\Gamma_0}}$ . \\

For checking well-definedness let also $\vartheta = \rund{\begin{array}{c|c} \sigma 1 & 0 \\ \hline
0 & F
\end{array}} \in G_0$~, $\sigma \in U(1)$ and $F \in U(r)$ , $\sigma^2 = \det F$ , such that $\vartheta \check \gamma \in \Gamma$ . Then again 
$\vartheta \eta^{- 1} \in \Gamma_0$ , which implies $\rund{\vartheta \eta^{- 1}}^{\abs{\Gamma_0}} = 1$ , and so 
$\sigma^{- \abs{\Gamma_0}} = \eps^{- \abs{\Gamma_0}}$ . $\Box$ \\

Let $\rho \in \{0, \dots, r\}$ be fixed. We intend to write the spaces $sM_k^\rho(\Gamma)$ and $sS_k^\rho(\Gamma)$ as spaces of global sections of 
holomorphic vector bundles on the compact Riemann surface 
$X = \left.\Gamma^\# \right\backslash H \cup \schweif{\text{cusps of } \left.\Gamma^\# \right\backslash H}$ .
For this purpose from now on given any $x \in \rz$ we will denote by $\floor{x} \in \zz$ and $\schweif{x} \in [0, 1[$ the unique numbers such that 
$x = \floor{x} + \schweif{x}$ . \\

We will construct holomorphic line bundles $L_k^0$ and $L_k^\twist$ and holomorphic vector bundles $M_k$ and $N_k$ on $X$ . Let us set

\begin{eqnarray*}
E_k^\rho &:=& \rund{T^* X}^{\otimes \floor{\frac{k + \rho}{2}} } \otimes L_k^0 \otimes L_k^\twist \otimes M_k  \phantom{123} \text{ and } \\
F_k^\rho &:=& \rund{T^* X}^{\otimes \floor{\frac{k + \rho}{2}} } \otimes L_k^0 \otimes L_k^\twist \otimes N_k \, .
\end{eqnarray*}

Then these bundles are supposed to have the following properties:

\begin{itemize}
\item[\{i\}] For all $U \subset X$ open $\Gamma^\hol\rund{U, E_\rho^k}$ ( $\Gamma^\hol\rund{U, F_\rho^k}$ ) is the space of super functions $f \in \O^\rho\rund{\pi_X^{- 1}(U)^{|r}}$ having $f|_\gamma = f$ for all $\gamma \in \Gamma$ and being bounded (vanishing) 
at each cusp of $\left.\Gamma^\# \right\backslash H$ belonging to $U$ , so in particular $sM_k^\rho(\Gamma) = H^0\rund{E_k^\rho}$ ( $sS_k^\rho(\Gamma) = H^0\rund{F_k^\rho}$ ),
\item[\{ii\}] $k \mapsto \deg L_k^\twist$ is a bounded function, and finally
\item[\{iii\}] the families $\rund{M_k}_{k \in \zz}$ and $\rund{N_k}_{k \in \zz}$ in fact consist of only finitely many vector bundles.
\end{itemize}

For this purpose let us split $k = l + 2 \abs{\Gamma_0} m$ , $l \in \schweif{0, \dots, 2 \abs{\Gamma_0} - 1}$ , $m \in \zz$~. We will construct all these bundles by local trivializations and identification of the fibres on the overlaps; $M_k$ and $N_k$ will have 
typical fibre $V_k^\rho$ . \\

{\it at regular points of $X$ }:

\begin{quote}
Let $z \in H$ be regular. Then its stabilizer in $\Gamma^\#$ is trivial and therefore the canonical projection $\pi_X: H \rightarrow X$ is locally biholomorphic at $z$ , so we take $\pi_X^{- 1}$ as a local chart of $X$ at $\overline z$ . Let 
$\gamma \in \Gamma$ . We have to indentify the fibres at the points $z$ and $\gamma^\# z$~: \\

for $L_k^0$ : $\rund{\gamma^\# z, s} \sim \rund{z, s}$ , \\

for $L_k^\twist$ : $\rund{\gamma^\# z, s} \sim \rund{z, \chi\rund{\check \gamma}^{2 m} s}$ , $s \in \cz$ , \\

and for $M_k$ and $N_k$ :

\[
\rund{\gamma^\# z, S} \sim \rund{z, j\rund{\check \gamma, z}^{2 \schweif{\frac{l + \rho}{2}} } \, S \ \varphi^\rho_l\rund{\check \gamma}} \, ,
\]

$S \in V_k^\rho$ , where $\check \gamma \in \check \Gamma$ is chosen to represent $\gamma^\#$ . Since

\[
S|_{\gamma} = j\rund{\check \gamma, z}^{k + \rho} \, S \ \varphi^\rho_k\rund{\check \gamma} = \chi\rund{\check \gamma}^{2 m} j\rund{\check \gamma, z}^{k + \rho} \, S \ \varphi^\rho_l\rund{\check \gamma}
\]

for all $S \in V_k^\rho$ , and since for $T^* X$ we have to identify

\[
\rund{\gamma^\# z, s} \sim \rund{z, \rund{\gamma^\#}'(z) \ s} \, ,
\]

where $\rund{\gamma^\#}'(z) = j\rund{\check \gamma, z}^2$ , we see that indeed $E_k^\rho$ and $F_k^\rho$ in \{i\} are given by the identification $\rund{\gamma^\# z, S} \sim \rund{z, S|_{\gamma}(z)}$ , $S \in V_k^\rho$ , and so \{i\} is true at regular points. \\
\end{quote}

{\it at elliptic points of $X$ }:

\begin{quote}
Let $z_0 \in H$ be elliptic. Then its stabilizer $\rund{\Gamma^\#}^{z_0}$ is cyclic of finite order $n \geq 2$ . Since the action of $\Aut H$ on $H$ is proper there exixts a $\rund{\Gamma^\#}^{z_0}$-invariant open neighbourhood $U \subset H$ of $z_0$ such 
that the canonical projection $\pi_{z_0}: \left.\rund{\Gamma^\#}^{z_0} \right\backslash H \rightarrow X$ restricted to $\left.\rund{\Gamma^\#}^{z_0} \right\backslash U$ is biholomorphic. Now let $c \in SL(2, \cz)$ be a Cayley transform mapping the unit 
disc $B \subset \cz$ biholomorphically onto $H$ such that $c \, 0 = z_0$ . We take

\[
\pi_X(U) \mathop{\longrightarrow}\limits^{\pi_{z_0}^{- 1}} \left.\rund{\Gamma^\#}^{z_0} \right\backslash U \mathop{\longrightarrow}\limits^{\overline c^{- 1}} \left.\spitz{e^{\frac{2 \pi i}{n}} } \right\backslash B \mathop{\longrightarrow}\limits^{w \mapsto w^n} B
\]

as a local chart of $X$ at $\overline{z_0}$ . Let $\rund{\Gamma^\#}^{z_0}$ be generated by $\gamma^\#$ , $\gamma \in \Gamma$ , such that $c^{- 1} \circ \gamma^\# \circ c \in \Aut(B)^0$ is the multiplication with $e^{\frac{2 \pi i}{n}}$ , and let again 
$\check \gamma \in \check \Gamma$ represent $\gamma^\#$ . \\

Let $e_\nu^l$ , $\nu = 1, \dots, \dim V_l^\rho$ , form a basis of $V_l^\rho = V_k^\rho$ , $e_\nu^l$ being an eigenvector of $j\rund{\check \gamma, z_0}^{l + \rho} \varphi^\rho_l\rund{\check \gamma}$ to the eigenvalue $e^{- 2 \pi i \sigma_\nu^l}$ , 
$\sigma_\nu^l \in \frac{1}{n} \zz$ , $\nu = 1, \dots, \dim V_l^\rho$ . Let $\chi\rund{\check \gamma}^2 = e^{- 2 \pi i \delta}$ , $\delta \in \frac{1}{n} \zz$ . Then every $e_\nu^l$ is at the same time an eigenvector of 
$j\rund{\check \gamma, z_0}^{k + \rho} \varphi^\rho_k\rund{\check \gamma}$ to the eigenvalue $e^{- 2 \pi i \rund{\sigma_\nu^l - \frac{m}{n} \abs{\Gamma_0} + m \delta}}$ . For defining the bundles we have to identify the fibres at the points 
$z \in U$ and $\rund{c^{- 1} z}^n \in B$ : \\

for $L_k^0$ : $\rund{z, s} \sim \rund{\rund{c^{- 1} z}^n, \rund{c^{- 1} z}^{n \, \floor{\frac{(k + \rho)(n - 1)}{2 n}} } s}$ , \\

for $L_k^\twist$ : $\rund{z, s} \sim \rund{\rund{c^{- 1} z}^n, \rund{c^{- 1} z}^{- n \schweif{m \delta}} s}$ , \\

and for $M_k$ and $N_k$ :

\[
\rund{z, e_\nu^l} \sim \rund{\rund{c^{- 1} z}^n, j\rund{c^{- 1}, z}^{- 2 \schweif{\frac{k + \rho}{2}} } \, \rund{c^{- 1} z}^{\Omega_\nu^k} e_\nu^l} \, ,
\]

where

\begin{eqnarray*}
\Omega_\nu^k &:=& n \rund{\schweif{m \delta} - \schweif{\sigma_\nu^l - \frac{m}{n} \abs{\Gamma_0} + m \delta}} + \floor{\frac{k + \rho}{2}}(n - 1) \\
&& \phantom{12} - n \floor{\frac{(k + \rho)(n - 1)}{2 n}} \, .
\end{eqnarray*}

Observe that $j\rund{c^{- 1} \check \gamma c, w} = j\rund{\check \gamma, z_0}$ for all $w \in B$ . Since for $T^* X$ we have to identify

\[
\rund{z, s} \sim \rund{\rund{c^{- 1} z}^n, \rund{c^{- 1}}'(z)^{- 1} \rund{c^{- 1} z}^{- n + 1} s}
\]

and $\rund{c^{- 1}}'(z) = j(c^{- 1}, z)^2$ we see that $E_k^\rho$ and $F_k^\rho$ in \{i\} are obtained by the identification

\[
\rund{z, e_\nu^l} \sim \rund{\rund{c^{- 1} z}^n, j\rund{c^{- 1}, z}^{- k - \rho} \, \rund{c^{- 1} z}^{- n \schweif{\sigma_\nu^l - \frac{m}{n} \abs{\Gamma_0} + m \delta}} e_\nu^l} \, .
\]

Now let $V \subset U$ be an open $\gamma^\#$-invariant neighbourhood of $z_0$ and $f \in \O(V)$ such that $f e_\nu^l$ is invariant under $|_\gamma$ . Then \\
$h := f\rund{c w} j\rund{c, w}^{k + \rho}$ fulfills

\[
h = h\rund{e^{\frac{2 \pi i}{n}} w} e^{- 2 \pi i \rund{\sigma_\nu^l - \frac{m}{n} \abs{\Gamma_0} + m \delta}} \, ,
\]

and so $\ord_0 h \geq n \schweif{\sigma_\nu^l - \frac{m}{n} \abs{\Gamma_0} + m \delta}$ . Therefore $h \, w^{- n \schweif{\sigma_\nu^l - \frac{m}{n} \abs{\Gamma_0} + m \delta}}$ is invariant under $w \mapsto e^{\frac{2 \pi i}{n}} w$ and still holomorphic at 
$w = 0$ , so $f e_\nu^l \in \Gamma^\hol\rund{\pi(V), E_k^\rho} = \Gamma^\hol\rund{\pi(V), F_k^\rho}$ . This shows \{i\} at elliptic points.
\end{quote}

{\it at cusps of $X$ }:

\begin{quote}
Let $z_0 \in \partial_{\pz^1} H$ be a cusp of $\left.\Gamma^\# \right\backslash H$ and $N^{z_0} \sqsubset \Aut H$ its associated nilpotent subgroup. Then $N^{z_0} \cap \Gamma^\#$ is infinite cyclic. Let $\gamma \in \Gamma$ and $g \in G$ such that
$\gamma^\#$ generates $N^{z_0} \cap \Gamma^\#$ and \\
$g_0 := g^{- 1} \gamma g$ is of the form (\ref{standardg_0}). Again choose an open $\gamma^\#$-invariant neighbourhood $U \subset H$ of $z_0$ such that  the canonical projection 
$\pi_{z_0} : \left.\spitz{\gamma^\#} \right\backslash H \cup \schweif{z_0} \rightarrow X$ restricted to $\left.\spitz{\gamma^\#} \right\backslash U \cup~\schweif{z_0}$ is biholomorphic. So

\begin{eqnarray*}
\pi_X(U) \cup \schweif{\overline{z_0}} \mathop{\longrightarrow}\limits^{\pi_{z_0}^{- 1}} \left.\spitz{\gamma^\#} \right\backslash U \cup \schweif{z_0}\mathop{\longrightarrow}\limits^{\overline{g^\#}^{- 1}} \spitz{z \mapsto z + 1} \backslash H \cup \{i \infty\} && \\
\mathop{\longrightarrow}\limits^{z \mapsto e^{2 \pi i z}} B \phantom{1234567890123456789012345678901234567890123456} &&
\end{eqnarray*}

is a local chart of $X$ at $\overline{z_0}$ . Let $\check g \in SL(2, \rz)$ represent $g^\#$ , and let $\check \gamma \in \check \Gamma$ represent $\gamma^\#$ such that $\check g^{- 1} \check \gamma \check g = \rund{\begin{array}{cc}
1 & 1 \\
0 & 1
\end{array}}$ . \\

Let again $e_\nu^l$ , $\nu = 1, \dots, \dim V_l^\rho$ , form a basis of $V_l^\rho = V_k^\rho$~, $e_\nu^l$ being an eigenvector of $\varphi^\rho_l\rund{\check \gamma}$ to the eigenvalue $e^{- 2 \pi i \sigma_\nu^l}$~, $\sigma_\nu^l \in \rz$ , 
$\nu = 1, \dots, \dim V_l^\rho$ . Let $\chi\rund{\check \gamma}^2 = e^{- 2 \pi i \delta}$ , $\delta \in \rz$ . Then $e_\nu^l$ is at the same time an eigenvector of $\varphi^\rho_k\rund{\check \gamma}$ to the eigenvalue $e^{- 2 \pi i \rund{\sigma_\nu^l + m \delta}}$ . We 
have to identify the fibres at the points $z \in U$ and $e^{2 \pi i {g^\#}^{- 1} z} \in B$ : \\

for $L_k^0$ : $\rund{z, s} \sim \rund{e^{2 \pi i {g^\#}^{- 1} z}, e^{2 \pi i \, {g^\#}^{- 1} z \, \floor{\frac{k + \rho}{2}} } s}$ , \\

for $L_k^\twist$ : $\rund{z, s} \sim \rund{e^{2 \pi i {g^\#}^{- 1} z}, e^{- 2 \pi i \, \schweif{m \delta} \, {g^\#}^{- 1} z} s}$ , \\

and finally for $M_k$ and $N_k$ :

\[
\rund{z, e_\nu^l} \sim \rund{e^{2 \pi i {g^\#}^{- 1} z}, j\rund{\check g^{- 1}, z}^{- 2 \schweif{\frac{k + \rho}{2}} } \, e^{2 \pi i \, \Omega_\nu^k \, {g^\#}^{- 1} z} \, e_\nu^l} \, ,
\]

where $\Omega_\nu^k := \schweif{m \delta} - \schweif{\sigma_\nu^l + m \delta}$ for $M_k$ and \\
$\Omega_\nu^k := \schweif{m \delta} - \rund{1 - \schweif{- \sigma_\nu^l - m \delta}}$ for $N_k$ . Observe that for all $x \in \rz$

\[
1 - \schweif{- x} = \left\{\begin{array}{ll}
\schweif{x} & \text{ if } x \notin \zz \\
1 & \text{ if } x \in \zz
\end{array} \right. \, .
\]

Again for proving \{i\} observe that $T^* X$ is given by the identification

\[
\rund{z, s} \sim \rund{e^{2 \pi i {g^\#}^{- 1} z}, \rund{{g^\#}^{- 1}}'(z)^{- 1} \, e^{- 2 \pi i {g^\#}^{- 1} z} \, s}
\]

and $\rund{{g^\#}^{- 1}}'(z) = j\rund{\check g^{- 1}, z}^2$ . So we obtain $E_k^\rho$ and $F_k^\rho$ in \{i\} by the identification

\[
\rund{z, e_\nu^l} \sim \rund{e^{2 \pi i {g^\#}^{- 1} z}, j\rund{\check g^{- 1}, z}^{- k - \rho} \, e^{- 2 \pi i \, \schweif{\sigma_\nu^l + m \delta} \, {g^\#}^{- 1} z} \, e_\nu^l}
\]

resp.

\begin{eqnarray*}
\rund{z, e_\nu^l} \phantom{123456789012345678901234567890123456789012345678} && \\
\sim \rund{e^{2 \pi i {g^\#}^{- 1} z}, j\rund{\check g^{- 1}, z}^{- k - \rho} \, e^{- 2 \pi i \, \rund{1 - \schweif{- \sigma_\nu^l - m \delta}} \, {g^\#}^{- 1} z} \, e_\nu^l} \, . &&
\end{eqnarray*}

Let again $V \subset U$ be an open $\gamma^\#$-invariant neighbourhood of $z_0$ and $f \in \O(V)$ such that $f e_\nu^l$ is invariant under $|_\gamma$ and bounded (vanishing) at $z_0$ . Then 
$h := f\rund{g^\# z} j\rund{\check g, z}^{k + \rho}$ is quasiperiodic $h = h\rund{z + 1} e^{- 2 \pi i \rund{\sigma_\nu^l + m \delta}}$ and bounded (vanishing) for $\Im z \leadsto + \infty$ . So $h \, e^{- 2 \pi i \schweif{\sigma_\nu^l + m \delta} z}$ (resp. 
$h \, e^{- 2 \pi i \rund{1 - \schweif{- \sigma_\nu^l - m \delta}} z}$ ) is invariant under $z \mapsto z + 1$ and bounded for $\Im z \leadsto + \infty$ . So $f e_\nu^l \in \Gamma^\hol\rund{\pi(V), E_k^\rho}$ (resp. $f e_\nu^l \in \Gamma^\hol\rund{\pi(V), F_k^\rho}$ ). \\
\end{quote}

Let $R$ denote the number of the elliptic points of $\left.\Gamma^\# \right\backslash H$ , and let $n_i$ be the period of the $i$th elliptic point. \\

For proving \{ii\} impose a metric on the holomorphic line bundle $L_k^\twist$ whose curvature is concentrated in small pairwise disjoint neighbourhoods of the elliptic points resp. cusps of $\left.\Gamma^\# \right\backslash H$ . It turns out that the total 
curvature of such a metric is bounded by $\pi \rund{R + S}$ . But the total curvature of any metric on the holomorphic line bundle $L_k^\twist$ is given by $\frac{\pi}{2} \deg L_k^\twist$ , so $\abs{\deg L_k^\twist} \leq 2 (R+ S)$ . \\

Now we prove \{iii\} :

\begin{quote}
Obviously there are at most $2 \abs{\Gamma_0}$ possibilities of how to indentify the fibres $V_k^\rho$ at $\gamma^\# z$ and $z$ , $z \in H$ regular. In the identification at an elliptic point $z_0$ since $\Omega_\nu^l \in \zz$ is of absolute value $< 2 n$ 
there are at most $4 n - 1$ possible values for $\Omega_\nu^l$ . In the identification at a cusp $z_0$ for fixed $l$ and $\nu$ there are at most $4$ possible values for $\Omega_\nu^l$ since on one hand $\abs{\Omega_\nu^l} < 2$ and on the other hand 
$\Omega_\nu^l \equiv - \sigma_\nu^l \, \mod \zz$ . $\Box$
\end{quote}

Obviously

\begin{eqnarray*}
\deg L_k^0 &=& \sum_{i = 1}^R \floor{\frac{(k + \rho) (n_i - 1)}{2 n_i }} + S \floor{\frac{k + \rho}{2}} \\
&=& \frac{k}{2} \rund{\sum_{i = 1}^R \rund{1 - \frac{1}{n_i}} + S} + O(1)
\end{eqnarray*}

for $k \leadsto \pm \infty$ . Let $g^*$ denote the genus of $X$ , so $\deg T^* X = 2 (g^* - 1)$ . A standard calculation using the total curvature of $X$ and the fact that $H$ is of constant curvature $- 1$ shows that 

\[
2 (g^* - 1) + \sum_{i = 1}^R \rund{1 - \frac{1}{n_i}} + S = \frac{\vol\rund{\left.\Gamma^\# \right\backslash H}}{2 \pi} > 0 \, .
\]

So we obtain the asymptotic behaviour

\begin{equation} \label{degree}
\deg\rund{\rund{T^* X}^{\otimes \floor{\frac{k + \rho}{2}} } \otimes L_k^0 \otimes L_k^\twist} = \frac{k}{2} \, \frac{\vol\rund{\left.\Gamma^\# \right\backslash H}}{2 \pi} + O(1) \leadsto \pm \infty
\end{equation}

for $k \leadsto \pm \infty$ .

\begin{lemma} \label{H1vanish} There exist $k_0, k_2 \in \zz$ such that 
\item[(i)] $H^1\rund{\rund{T X \otimes E_k^\rho}^*} = 0$ for all $k \leq k_0$ , and
\item[(ii)] for all $k \geq k_2$ : $H^1\rund{E_k^\rho} , H^1\rund{F_k^\rho} = 0$ , and $\Gamma^\hol\rund{\diamondsuit, E_k^\rho}$ is generated by global sections.
\end{lemma}

{\it Proof:} By \{iii\} we may assume that $M_k$ and $N_k$ are independent of $k$ . So we obtain the result combining (\ref{degree}) and lemma 7.1 b) of \cite{Hart} , which says that given any coherent sheaf $F$ on a non-singular projective curve $X$ , 
there is an integer $d_0$ such that if $L$ is a line bundle over $X$ of degree $\geq d_0$ , then $F \otimes L$ is generated by global sections, and $H^1(F \otimes L) = 0$ . $\Box$ \\

{\it Now we prove theorem \ref{main}}: (i) Serre duality tells us that

\[
sM_k(\Gamma) = H^0\rund{E_k^\rho} \simeq H^1\rund{\rund{T X \otimes E_k^\rho}^*}^* \, ,
\]

which is $0$ if $k \leq k_0$ , $k_0 \in \zz$ be given by lemma \ref{H1vanish}.

(ii) By the Riemann Roch theorem applied to $E_k^\rho$ , which is of rank \\
$n_k := \dim V_k^\rho$ , we obtain

\[
\dim H^0\rund{E_k^\rho} - \dim H^1\rund{E_k^\rho} = c_1\rund{E_k^\rho} - n_k \rund{g^* - 1} \, ,
\]

where $c_1\rund{E_k^\rho} = \deg \bigwedge^{n_k} E_k^\rho$ denotes the first Chern class of $E_k^\rho$ . But $\dim H^1\rund{E_k^\rho} = 0$ for $k \geq k_2$ , $k_2 \in \zz$ be given by lemma \ref{H1vanish}, and

\begin{eqnarray*}
c_1\rund{E_k^\rho} &=& n_k \, \deg \rund{\rund{T^* X}^{\otimes \floor{\frac{k + \rho}{2}} } \otimes L_k^0 \otimes L_k^\twist} + \deg \bigwedge\nolimits^{n_k} M_k \\
&=& n_k \rund{\frac{k}{2} \phantom{1} \frac{\vol\rund{\left.\Gamma^\# \right\backslash H}}{2 \pi} + O(1)} \, ,
\end{eqnarray*}

which gives the asymptotic formula. For proving the inequality take a cusp $z_0 \in \partial_{\pz^1} H$ and the associated basis $e_\nu^\rho$ , 
$\nu = 1, \dots, \dim V_l^\rho$ , of $V_l^\rho$~. For $\nu = 1, \dots, \dim V_l^\rho$ associate the coefficient of $e_\nu^\rho$ in 
$f\rund{\overline{z_0}} \in V_l^\rho$ if $\sigma_\nu^l + m \delta \in \zz$ and $0$ otherwise to every $f \in sM_k^\rho(\Gamma) = H^0\rund{E_k}$. Putting all 
cusps together yields a linear map $sM_k(\Gamma) \rightarrow \cz^{\ S \dim V_l^\rho}$ with kernel $sS_k(\Gamma)$ .

(iii) Let $k_2 \in \zz$ be given by lemma \ref{H1vanish} with respect to $\rho = 1$ and let $k_1 \in \zz_{\geq k_2}$ be given such that for all $k \geq k_1$ the holomorphic line bundle $E_k^0$ is already very ample. Let $k \geq k_1$ be arbitrary. Then of 
course

\[
\Phi^\#: \left.\Gamma^\# \right\backslash H \hookrightarrow X \rightarrow \pz^m
\]

is already an embedding. Now let $z_0 \in H$ be arbitrary. Without loss of generality we may assume that $f_0\rund{z_0} \not= 0$ . So using the 0th standard local super chart of $\pz^{m|n}$ , $\Phi$ is given by the tuple 
$\frac{1}{f_0} \rund{f_1, \dots, f_m, \lambda_1, \dots, \lambda_n}$ in some neighbourhood of $z_0$ . Since $\Gamma^\hol\rund{\diamondsuit, E_k^1}$ is generated by global sections according to lemma \ref{H1vanish} (ii) we see that 
$\rund{sD \ \Phi}^\#\rund{z_0}$ is injective, and so $\Phi$ is a super embedding by the super inversion theorem.

If $\left.\Gamma^{\#} \right\backslash H$ has no cusps then it is compact, and so algebraicity follows from a super version of Chow's theorem, see theorem 6 of \cite{LeBrun} . $\Box$ \\

\section{$\mathcal{P}$-lattices}

Let $\P = \P_0 \oplus \P_1$ be a real finite dimensional unital associative and super commutative super algebra having a unique maximal ideal $\m$ (so $\P$ is local, and automatically $\P_1 \sqsubset \m$ and $\m$ is graded), $\m^N = 0$ for some 
$N \in \nz$ , and a canonical projection ${}^{\#'} : \P \rightarrow \P / \m \simeq \rz$ . Examples are $\P = \bigwedge\rund{\rz^{N - 1}}$ and $\P = \rz[X] \left/ \rund{X^N}\right.$~, the second being purely even. As promised, for a lattice 
$\Gamma \sqsubset G = \G^\#$ we will now discuss super deformations of the embedding $\Gamma \hookrightarrow \G$ 'parametrized' by the generators of $\P$ . We will call such super deformations $\P$-lattices and give a precise definition in a moment.

\begin{defin} [$\P$-points] \end{defin}

Let $\M = (M, \S)$ be a real super manifold of super dimension $(m, n)$ , $M = \M^\#$ being an ordinary smooth $n$-dimensional manifold and $\S$ a sheaf of unital associative super commutative super algebras on $M$ , locally 
$\simeq \C^\infty \otimes \bigwedge\rund{\rz^n}$ . Then a $\P$-point of $\M$ is a morphism $A$ of from $(\{0\}, \P)$ to $\M$ as  ringed spaces. Here an equivalent definition: A pair $A := (a, \a)$ where $a \in M$ is an ordinary point and 
$\a: \S_a \rightarrow \P$ , where $\S_a$ denotes the stalk of $\S$ at $a$ , is called a $\P$-point of $\M$ . $A^{\#'} := a \in M$ is called the relative body of $A$ . We write $A \in_\P \M$ . The set of $\P$-points of $\M$ is denoted by $\M^\P$ . \\

Having chosen local super coordinates on $\M$ , the $\P$-points of $\M$ lying in the range of these are in 1-1-correspondence with tuples $\rund{a_1, \dots, a_m, \alpha_1, \dots, \alpha_n} \in \P_0^{\oplus m} \oplus \P_1^{\oplus n}$ , and in this notation the 
relative body is given by $\rund{a_1^{\#'}, \dots, a_m^{\#'}} \in \rz^{\oplus m}$ . If $N = 2$ (infinitesimal super deformation) we have a 1-1-correspondence between $\M^\P$ and pairs $(a, v)$ where $a \in M$ and $v \in \rund{sT_a \M \otimes \P}_0$ , 
$sT_a \M$ denoting the super tangent space of $\M$ at $a$ . \\

Obviously every super morphism between the real super manifolds \\
$\M = (M, \S)$ and $\N = (N, \T)$ induces a map $\M^\P \rightarrow \N^\P$ . So we obtain a whole functor from the category of real super manifolds to the category of sets, 
and this functor restricts to a functor from the category of real super Lie groups to the category of groups. Indeed, given a real super Lie group $\G$ with body $G$ , the multiplication super morphism $m: \G \times \G \rightarrow \G$ turns the set $G^\P$ 
of all $\P$-points of $G$ into a group via $g h := m(g, h)$ for all $g, h \in_\P G$~, and clearly ${}^{\#'}: \G^\P \rightarrow G \, , \, g \mapsto g^{\#'}$ is a group epimorphism.

\begin{defin} [ $\P$-lattices] \end{defin}

Let $\G$ be a real super Lie group with body $G$ and $\Upsilon \sqsubset \G^\P$ be a subgroup. $\Upsilon$ is called a $\P$-lattice of $\G$ iff

\begin{itemize}
\item[\{i\}] $\Upsilon^{\#'} := \schweif{\left.\gamma^{\#'} \, \right| \, \gamma \in \Upsilon} \sqsubset G$ is an ordinary lattice, called the relative body of $\Upsilon$ , and
\item[\{ii\}] ${}^{\#'} : \Upsilon \rightarrow \Upsilon^{\#'} \, , \, \gamma \mapsto \gamma^{\#'}$ is bijective and so automatically an isomorphism.
\end{itemize}
 
Of course given a $\P$-lattice $\Upsilon$ of $\G$ with relative body $\Gamma \sqsubset G$ and $g \in_\P \G$ with $g^{\#'} = 1$ we get another $\P$-lattice $g \Upsilon g^{- 1}$ of $\G$ with same relative body $\Gamma$ , and we are interested in classifying 
all the conjugacy classes for given $\Gamma$~. If $N = 2$ they are in 1-1-correspondence with $\rund{H^1\rund{\Gamma, \g} \otimes \m}_0$ , $\Gamma$ acting on the super Lie algebra $\g$ of $\G$ by $s\Ad$ , compare with the classical case for example in 
\cite{Ragh}. \\

One is also interested in the question if it is always possible to extend a given local super deformation of $\Gamma$ to higher degree $N$ of nilpotency: Let $\Q := \P \left/ \m^{N - 1}\right.$ . Then $\Q$ fulfills the same properties as $\P$ with maximal ideal 
$\n := \m \left/ \m^{N - 1}\right.$ , $\n^{N - 1} = 0$ . The canonical projection ${}^\natural: \P \rightarrow \Q$ obviously induces a map respecting ${}^{\#'}$ from $\P$-points of a super manifold $\M$ to its $\Q$-points. Now given a $\Q$-lattice $\Upsilon$ of a 
super Lie group $\G$ , does there exist a $\P$-lattice $\hat \Upsilon$ such that ${\hat \Upsilon}^\natural = \Upsilon$ ? As in the classical case the answer is yes if $H^2(\Gamma, \g) = 0$ , and the converse is false. \\

Given a $\P$-lattice $\Upsilon$ of a super Lie group $\G$ and $\gamma \in \Upsilon$ such that $\rund{\gamma^{\#'}}^n = 1$ for some $n \in \nz \setminus \{0\}$ , automatically $\gamma^n = 1$ , and so by the following lemma $\gamma$ is conjugate to 
$\gamma^{\#'}$ : there exists $g \in_\P \G$ such that $g^{\#'} = 1$ and $\gamma = g \gamma^{\#'} g^{- 1}$ .

\begin{lemma} \label{finord} Let $\G$ be a super Lie group with body $G$ and super Lie algebra $\g$ , and let $n \in \nz \setminus \{0\}$ . Then the equation $g^n = 1$ defines sub super manifolds $\M$ of $\G$ whose bodies are precisely the connected 
components of $M := \schweif{g \in G \, \left| \, g^n = 1\right.}$ . Let $g_0 \in M$ and $V \sqsubset \g$ be a graded complement of $\z_\g\rund{g_0}$ . Then $\exp(\chi) \, g_0 \, \exp(- \chi)$ locally at $0 \mapsto g_0$ defines a super diffeomorphism 
$V \rightarrow \M$ , $V$ regarded as a real super manifold.
\end{lemma}

{\it Proof:} Let the super morphism $\Omega: \g \rightarrow \G$ be given by $\exp\rund{\chi_V} \, g_0 \exp\rund{\chi_\z} \, \exp\rund{- \chi_V}$ , where $\chi_V$ and $\chi_\z$ denote the projections on $\g$ along the splitting $\g = V \oplus \z_\g\rund{g_0}$ , 
$\g$ treated as a real super manifold. Then a straight forward calculation shows that the super differential $sD \Omega(0)$ is bijective, and so $\Omega$ is a super diffeomorphism locally at $0$ . Now $\Omega^n$ defines a super morphism 
$\Psi: \g \rightarrow \G$ having $\left.\Psi\right|_V \equiv 1$ and $sD \, \Psi (0) |_{\z_\g\rund{g_0}}$ is injective. Therefore the equation $\Psi = 1$ locally at $0$ defines the sub super manifold $V$ of $\g$ . $\Box$ \\

From now on let again $\G$ be the real sub super Lie group of $\rund{GL(2, \cz) \times GL(r, \cz)}^{| 4 r}$ from section 1 given by the equations $g I g^* = I$ and $\Ber g = 1$ . Then $\G^\P$ is the set of all even super matrices $\rund{\P^{(2|r) \times (2|r)} }_0$ 
(even entries in the diagonal, odd entries in the off-diagonal blocks) fulfilling these two equations, and the product of two of them can be computed via ordinary matrix multiplication. Of course the action $\alpha: \G \times H^{|r} \rightarrow H^{|r}$ induces a 
group homomorphism from $\G^\P$ into the group of $\P$- super automorphisms of $H^{|r}$ respecting ${}^{\#'}$ :

\begin{defin} [ $\P$- super automorphisms of $H^{|r}$ ] \end{defin}

An automorphism $\Phi$ of the ringed space $\rund{H, \P^\cz \boxtimes \rund{\O_H \otimes \bigwedge \rund{\cz^r}} }$ , \\
$\P^\cz \boxtimes \rund{\O_H \otimes \bigwedge \rund{\cz^r}}$ treated as a sheaf of unital $\zz_2$-graded $\P^\cz$-modules, is called a $\P$- super automorphism of $H^{|r}$ . Clearly the projection \\
${}^{\#'} : \P \rightarrow \rz$ induces an embedding $H^{|r} \hookrightarrow \rund{H, \P^\cz \boxtimes \rund{\O_H \otimes \bigwedge \rund{\cz^r}} }$ as ringed spaces whose underlying map 
$H \rightarrow H$ is the identity. The unique super automorphism $\Phi^{\#'}$ of $H^{|r}$ such that

\[
\begin{array}{ccc}
\phantom{12345} H^{|r} & \hookrightarrow & \rund{H, \P^\cz \boxtimes \rund{\O_H \otimes \bigwedge \rund{\cz^r}} } \\
\Phi^{\#'} \downarrow & \circlearrowleft & \downarrow \Phi \\
\phantom{12345} H^{|r} & \hookrightarrow & \rund{H, \P^\cz \boxtimes \rund{\O_H \otimes \bigwedge \rund{\cz^r}} }
\end{array}
\]

is called the relative body of $\Phi$ . \\

In practice $\P$- super automorphisms of $H^{|r}$ are given by tuples $\rund{f, \lambda_1, \dots, \lambda_r} \in \rund{\P^\cz \boxtimes \O\rund{H^{|r}} }_0 \oplus \rund{\P^\cz \boxtimes \O\rund{H^{|r}} }_1^{\oplus r}$~, and in this notation the relative body is 
given by the tuple \\
$\rund{f^{\#'}, \lambda_1^{\#'}, \dots, \lambda_r^{\#'}} \in \O\rund{H^{|r}}_0 \oplus \O\rund{H^{|r}}_1^{\oplus r}$ , where ${}^{\#'}$ denotes the complexification and right- $\O_H \otimes \bigwedge\rund{\cz^r}$ -linear extension of the projection 
${}^{\#'}: \P \rightarrow \rz$ . \\

Do not mix up the body ${}^\#$ and the relative body ${}^{\#'}$ : \\

Given some $g \in_\P \G$ and some $\P$- super automorphism $\Phi$ of $H^{|r}$ , $g^{\#'}$ , the {\bf relative body} of $g$ , is an ordinary point of $G$ , while the {\bf body} $g^\#$ of $g$ by definition coincides with 
the body of $g^{\#'}$ and is an element of $\Aut H$ .

The {\bf relative body} $\Phi^{\#'}$ of $\Phi$ is {\bf still} a {\bf super} automorphism of $H^{|r}$ , while the {\bf body} $\Phi^\#$ of $\Phi$ is the underlying {\bf ordinary} automorphism of $H$ . $g$ and $\Phi$ 
are local super deformations over $\P$ resp. $\P^\cz$ of their {\bf relative} bodies $g^{\#'}$ and $\Phi^{\#'}$ .

Taking the {\bf relative body} ${}^{\#'}$ means set all generators of $\P$ to zero, taking the {\bf body} ${}^\#$ means set everything to zero which is nilpotent, the generators of $\P$ {\bf and} the odd coordinates $\zeta_1, \dots, \zeta_r$ on $H^{|r}$ . \\

In the end of this section let us discuss - for the lattices $\Gamma \sqsubset G$ of the examples \ref{exordlat} - $H^1(\Gamma, \g)$ and the $\P$-lattices $\Upsilon$ of $\G$ with relative body $\Gamma$ : \\

First we observe that after identification $\g_1 \simeq \cz^{r \times 2}$

\begin{eqnarray*}
\z_{\g_0}\rund{\gamma_0} &=& \sl(2, \rz) \oplus \z_{\u(r)}\rund{E_0} \, , \\
\z_{\g_1}\rund{\gamma_0} &=& \Eig_{\eps_0}\rund{E_0}^{\oplus 2} \, .
\end{eqnarray*}

\begin{itemize}
\item[$\spitz{\text{i}}$] By lemma \ref{finord} we see that if $V$ is a graded complement of \\
$\z_\g\rund{\gamma_0, \hat R} + \z_\g\rund{\gamma_0, \hat S}$ in $\z_\g\rund{\gamma_0}$ then the conjugacy classes of $\P$-lattices $\Upsilon$ of $\G$ with relative body $\Gamma$ are in 1-1-correspondence with $\rund{V \otimes \m}_0$ via the 
assignment $\chi \mapsto \spitz{\gamma_0, \hat R, \exp(\chi) \hat S \exp(- \chi)}$ . So there are no obstructions for extending a local super deformation of $\Gamma$ to higher degree of nilpotency, and

\[
H^1(\Gamma, \g) \simeq \z_\g\rund{\gamma_0} \, \left/ \, \rund{\z_\g\rund{\gamma_0, \hat R} + \z_\g\rund{\gamma_0, \hat S}} \right. \, .
\]

\begin{eqnarray*}
H^1\rund{\Gamma, \g_0} &\simeq& \sl(2, \rz) \, \left/ \,  \rund{\z_{\sl(2, \rz)}(R) + \z_{\sl(2, \rz)}(S)} \right. \\
&& \phantom{12} \oplus \z_{\u(r)}\rund{E_0} \left/ \rund{\z_{\u(r)}\rund{E_0, E} + \z_{\u(r)}\rund{E_0, F}} \right. \, ,
\end{eqnarray*}

where the first summand is of dimension $1$ . Since $E$ and $F$ commute with $E_0$ we may define $\varphi, \psi \in GL\rund{\Eig_{\eps_0}\rund{E_0}^{\oplus 2}}$ as $u \mapsto E u R^{- 1}$ resp. $u \mapsto F u S^{- 1}$ . Then

\[
H^1\rund{\Gamma, \g_1} \simeq \left.\Eig_{\eps_0}\rund{E_0}^{\oplus 2} \right/ \rund{\Eig_{\eps}(\varphi) + \Eig_{\eta}(\psi)} \, ,
\]

which has maximal real dimension $4 r = \dim \g_1$ if for example \\
$E_0 = E = 1$ , $\det F = - 1$ , $F$ has no real eigenvalues, $\eps_0 = \eps = 1$ and $\eta = i$ .

\item[$\spitz{\text{ii}}$] {\it First case:} $m \geq 1$ . Again by lemma \ref{finord} it is enough to consider $\P$-lattices $\Upsilon$ of $\G$ having $\gamma_0 \in \Upsilon$ modulo conjugation with $g \in \cZ_\G\rund{\gamma_0}$~, $g^{\#'} = 1$ , where 
$\cZ_\G\rund{\gamma_0}$ denotes the centralizer of $\gamma_0$ in $\G$ , which is a sub super Lie group of $\G$ with super Lie algebra $\z_\g\rund{\gamma_0}$ . Obviously these lattices are given by $\P$-points of 
$\cZ_\G\rund{\gamma_0}^{2 g^* + m - 1}$ with body $\rund{\hat A_k, \hat B_k, \hat C_l}_{k = 1, \dots, g^*, l = 1, \dots, m - 1}$ (from now on we will drop the index), but not in 1-1-correspondence, we still have to devide out conjugation. However we observe 
that again there are no obstructions for extending a local super deformation of $\Gamma$ to higher degree of nilpotency , and

\[
\sdim H^1(\Gamma, \g) = \rund{2 g^* + m - 2} \sdim \z_\g\rund{\gamma_0} \, + \, \sdim \z_\g\rund{\gamma_0, \hat A_k, \hat B_k, \hat C_l} \, .
\]

Since $\z_\g\rund{\gamma_0, \hat A_k, \hat B_k, \hat C_l} = \z_{\u(r)}\rund{E_0, E_k, F_k, C_l}$ is purely even we obtain

\begin{eqnarray*}
\dim H^1\rund{\Gamma, \g_0} &=& \rund{2 g^* + m - 2} \rund{3 + \dim \z_{\u(r)} \rund{E_0}} \\
&& \phantom{12} + \dim \z_{\u(r)}\rund{E_0, E_k, F_k, C_l} \, , \\
\dim H^1\rund{\Gamma, \g_1} &=& 2 \rund{2 g^* + m - 2} \dim \Eig_{\eps_0}\rund{E_0} \, .
\end{eqnarray*}

{\it Second case:} $m = 0$ . Then $\gamma_0 = 1$ , and with the super morphisms $\Phi: \G^{2 g^*} \rightarrow \G$ defined as $\eckig{g_1, h_1} \cdots \eckig{g_{g^*}, h_{g^*}}$ and \\
$\Psi: \G \rightarrow \G^{2 g^*}$ defined as $\rund{g \hat A_k g^{- 1}, g \hat B_k g^{- 1}}$ we see that \\
$H^1(\Gamma, \g) = \left.\ker sD \Phi\rund{\hat A_k, \hat B_k} \right/ \Im sD \Psi(1)$ , so

\begin{eqnarray*}
\sdim H^1(\Gamma, \g) &=& \rund{2 g^* - 1} \sdim \g - \sdim \Im sD \Phi\rund{\hat A_k, \hat B_k} \\
&& \phantom{12} + \sdim \z_\g\rund{\hat A_k, \hat B_k} \, .
\end{eqnarray*}

Some longer calculations show that

\[
\Im sD \Phi\rund{\hat A_k, \hat B_k}_0 = \sl(2, \rz) \oplus \z_{\su(r)}\rund{E_k, F_k}^\perp \, ,
\]

where ${}^\perp$ is taken with respect to the Killing form on $\su(r)$ , $\Im sD \Phi\rund{\hat A_k, \hat B_k}_1 = \g_1$ , and $\z_\g\rund{\hat A_k, \hat B_k} = \z_{\u(r)}\rund{E_k, F_k}$ is purely even. So in the end

\begin{eqnarray*}
\dim H^1\rund{\Gamma, \g_0} &=& 2 \rund{g^* - 1} \rund{3 + r^2} + 2 \dim \z_{\u(r)}\rund{E_k, F_k} \, , \\
\dim H^1\rund{\Gamma, \g_1} &=& 8 \rund{g^* - 1} r \, .
\end{eqnarray*}

In contrast to the case $m \geq 1$ here one can construct examples with obstructions for extending a local super deformation of $\Gamma$ to higher degree of nilpotency.

\end{itemize}

\section{local sheaf deformation}

Throughout this section let $X$ be a topological space and $\P$ a finite dimensional unital associative super algebra over a field $K$ having a unique maximal ideal $\m$~, 
$\m^N = 0$ for some $N \in \nz$ , and a canonical projection ${}^{\#'} : \P \rightarrow \P / \m \simeq K$ . Let $\E$ be a sheaf of left-$\P$-modules over $X$ such that locally

\[
\begin{array}{ccc}
\E & \simeq & \P \otimes E \\
{}_{{}^{\#'}} \searrow & \circlearrowleft & \swarrow_{{}^{\#'} \otimes \id} \\
 & E & 
\end{array} \, ,
\]

where $E := \E / \m \E$ , which is a sheaf of $K$-vectorspaces over $X$ , and \\
${}^{\#'} : \E \rightarrow E$ denotes the canonical projection. Then of course $\E$ can be given by an open cover $\rund{U_i}_{i \in I}$ of $X$ , isomorphisms 

\[
\begin{array}{ccc}
\E|_{U_i} & \simeq & \P \otimes E|_{U_i} \\
{}_{{}^{\#'}} \searrow & \circlearrowleft & \swarrow_{{}^{\#'} \otimes \id} \\
 & E|_{U_i} & 
\end{array}
\]

and transition functions

\[
\varphi_{i j} = \id + A_{i j} : \P \otimes E|_{U_i \cap U_j} \rightarrow \P \otimes E|_{U_i \cap U_j} \, ,
\]

$i, j \in I$ , between them, where $A_{i j} : \P \otimes E|_{U_i \cap U_j} \rightarrow \m \otimes E|_{U_i \cap U_j}$ are left-$\P$-linear maps. Obviously $\left.\m^{N - 1} E\right|_{U_i \cap U_j} \sqsubset \ker A_{i j}$ , so these local isomorphisms glue together to a 
canonical global isomorphism $\m^{N - 1} \E \simeq \m^{N - 1} \otimes E$~.

\begin{lemma} \label{globsect}
Let $d := \dim E(X) < \infty$ . Then

\[
d \leq \dim_K \E(X) \leq d \, \dim_K \P  \, ,
\]

and equivalent are

\begin{itemize}
\item[(i)] $\dim_K \E(X) = d \, \dim_K \P$ ,
\item[(ii)] there exist $f_1, \dots, f_d \in \E(X)$ such that $\rund{f_1^{\#'}, \dots, f_d^{\#'}}$ is a basis of $E(X)$ ,
\item[(iii)] $\E(X)$ is a free $\P$-module of rank $d$ .
\end{itemize}

Furthermore if (ii) is valid then $\rund{f_1, \dots, f_d}$ is a $\P$-basis of $\E(X)$ , and the assignment $f_\delta \mapsto f_\delta^{\#'}$ , $\delta = 1, \dots, d$ , induces a $\P$-module isomorphism

\[
\begin{array}{ccc}
\E(X) & \simeq & \P \otimes E(X) \\
{}^{\#'} \searrow & \circlearrowleft & \swarrow {}^{\#'} \otimes \id \\
 & E(X) & 
\end{array} \, .
\]

\end{lemma}

{\it Proof:} The first inequality is of course trivial if $\m = 0$ .  For $\m \not= 0$ let $N' \in \nz$ be maximal such that $\m^{N'} \not= 0$ . Then $\m^{N'} \otimes E(X) = \m^{N'} \E(X) \sqsubset \E(X)$ , which proves the first inequality. \\

The second inequality, the implication (i) $\Rightarrow$ (ii) and the last statement will be proven by induction on $N \in \nz \setminus \{0\}$ . If $N = 1$ then $\m = 0$ and all statements are trivial.

Now assume $\m^{N + 1} = 0$ . Then define $\Q := \P \left/ \m^N \right.$ , which has the unique maximal ideal $\n := \m \left/ \m^N \right.$ , $\n^N = 0$ and $\Q / \n \simeq K$ , and let ${}^\natural : \P \rightarrow \Q$ be the canonical projection. Let 
$\E^\natural := \E \left/ \m^N \E \right.$ and

\[
{}^\natural : \E(X) \rightarrow \E^\natural(X)
\]

be the linear map induced by the canonical sheaf projection $\E \rightarrow \E^\natural$ . Its kernel is $\m^N \E(X) = \m^N \otimes E(X)$ . By induction hypothesis \\
$\dim_K \E^\natural(X) \leq d \, \dim_K \Q$ , and so

\[
\dim_K \E(X) \leq d \, \dim_K \Q + d \, \dim_K \m^N = d \, \dim_K \P  \, ,
\]

which proves the second inequality.

For proving the implication (i) $\Rightarrow$ (ii) assume $\dim_K \E(X) = d \, \dim_K \P$~. Then since $\dim_K \P = \dim_K \Q + \dim_K \m^N$ , \\
$\dim_K \rund{\m^N \otimes E(X)} = d \, \dim_K \m^N$ and $\dim_K \E^\natural(X) \leq d \, \dim_K \Q$ we see that necessarily 

\[
{}^\natural : \E(X) \rightarrow \E^\natural(X)
\]

is surjective and $\dim_K \E^\natural(X) = d \, \dim_K \Q$ . So by induction hypothesis and surjectivity there exist $f_1, \dots, f_d \in \E(X)$ such that $\rund{f_1^{\#'} , \dots, f_d^{\#'}}$ is a basis of $E(X)$ , which proves (ii) .

For proving the last statement let $f_1, \dots, f_d \in \E(X)$ such that $\rund{f_1^{\#'} , \dots, f_d^{\#'}}$ is a basis of $E(X)$ . Then by induction hypothesis $\rund{f_1^\natural, \dots, f_d^\natural}$ is a $\Q$-basis of $\E^\natural(X)$ . For 
proving that $\rund{f_1, \dots, f_d}$ spans $\E(X)$ over $\P$ let $F \in \E(X)$ . Then there exist $a_1, \dots, a_d \in \P$ such that

\[
F^\natural = a_1^\natural f_1^\natural + \dots + a_d^\natural f_d^\natural  \, ,
\]

and so

\[
\Delta := F - a_1^\natural f_1^\natural - \dots - a_d^\natural f_d^\natural \in \m^N \E(X) = \m^N \otimes E(X) \, .
\]

Since $\rund{f_1^{\#'}, \dots, f_d^{\#'}}$ is a basis of $E(X)$ we see that there exist $b_1, \dots, b_d \in~\m^N$ such that

\[
\Delta = b_1 \otimes f_1^{\#'} + \dots + b_d \otimes f_d^{\#'} = b_1 f_1 + \dots + b_d f_d  \, ,
\]

and so

\[
F = \rund{a_1 + b_1} f_1 + \dots + \rund{a_d + b_d} f_d  \, .
\]

For proving linear independence let $a_1, \dots, a_d \in \P$ such that

\[
a_1 f_1 + \dots + a_d f_d = 0 \, .
\]

Then $a_1^\natural f_1^\natural + \dots + a_d^\natural f_d^\natural = 0$ in $\E^\natural(X)$ , and so $a_1^\natural = \dots = a_d^\natural = 0$ . Therefore $a_1, \dots, a_d \in \m^N$ , and this means

\[
0 = a_1 f_1 + \dots + a_d f_d = a_1 \otimes f_1^{\#'} + \dots + a_d \otimes f_d^{\#'} \, .
\]

Since $f_1^{\#'}, \dots, f_d^{\#'}$ are linearly independent over $K$ we get $a_1 = \dots = a_d = 0$ . \\

Now (ii) $\Rightarrow$ (iii) follows from the last statement, and (iii) $\Rightarrow$ (i) is of course trivial. $\Box$ \\

The crutial question is now: Given an element $f \in E(X)$ , is it possible to adapt $f$ to the local deformation $\E$ of $E$ , precisely, is it possible to construct $\widetilde f \in \E(X)$ such that $\widetilde f^{\#'} = f$ ?

\begin{lemma} \label{adapt}
Assume that $H^1(X, E) = 0$ . Then for all $f \in E(X)$ there exists $\widetilde f \in \E(X)$ such that ${\widetilde f}^{\#'} = f$ .
\end{lemma}

{\it Proof:} via induction on $N \in \nz \setminus \{0\}$ . If $N = 1$ again the statement is trivial.

Now assume $\m^{N + 1} = 0$ . Again define $\Q := \P \left/ \m^N \right.$ with unique maximal ideal $\n := \m \left/ \m^N \right.$ and canonical 
projection ${}^\natural : \P \rightarrow \Q$~. Let $f \in E(X)$~. Then by induction hypothesis there exists $\widetilde f' \in \E^\natural(X)$ such that 
$\widetilde f'^{\#'} = f$~. Since $\E^{\natural}$ is given by local isomorphisms $\E^\natural \simeq \Q \otimes E$ 
with transition functions $\varphi_{i j}^\natural = \id + A_{i j}^\natural : \Q \otimes E|_{U_i \cap U_j} \rightarrow \Q \otimes E|_{U_i \cap U_j}$ we see that $\widetilde f'$ is given by sections $f|_{U_i} - \sigma_i^\natural \in \Q \otimes E|_{U_i}$ , where 
$\sigma_i \in \m \otimes E\rund{U_i}$ , $i \in I$ . Using $\varphi_{i j}^\natural\rund{f - \sigma_i^\natural} =  f - \sigma_j^\natural$ on $U_i \cap U_j$ , an easy 
calculation shows that

\[
a_{i j} := \varphi_{i j}\rund{f - \sigma_i} - f + \sigma_j \in \m^N \otimes E\rund{U_i \cap U_j} \, ,
\]

$i, j \in I$ , define a cocycle in $\m^N \otimes Z^1\rund{\rund{U_i}_{i \in I} , E}$ . Since by assumption $H^1(X, E) = 0$ we see that after maybe some refinement of the open cover $\rund{U_i}_{i \in I}$ we may assume that there exist 
$\tau_i \in \m^N \otimes E\rund{U_i}$ , $i \in I$ , such that $a_{i j} = \tau_i - \tau_j$ . Again an easy calculation shows that $f - \sigma_i - \tau_i \in \P \otimes E\rund{U_i}$~, $i \in I$ , glue together to an element $ \widetilde f \in \E(X)$ having 
$\widetilde f^{\#'} = f$ . $\Box$ 

\section{Super automorphic forms for $\mathcal{P}$-lattices}

Let again $\P$ be as in section 3 and $k \in \zz$ . For $g \in_\P \G$ , $U \subset H$ open and $f \in \O\rund{\rund{g^\# U}^{|r}}$ there is little hope that $f\rund{g \rund{\Z}} j\rund{g, \Z}^k$ will lie in $\O\rund{U^{|r}}$ . However,

\[
|_{g, k} : \P^\cz \boxtimes \O\rund{\rund{g^\# U}^{|r}} \rightarrow \P^\cz \boxtimes \O\rund{U^{|r}} \, , \, f \mapsto f\rund{g \rund{\Z}} j\rund{g, \Z}^k
\]

defines a $\zz_2$-graded $\P^\cz$-linear map, and so in particular we obtain a right representation of $\G^\P$ on $\P^\cz \boxtimes \O\rund{H^{|r}}$ . \\

For defining super automorphic resp. cusp forms for a $\P$-lattice $\Upsilon$ of $\G$ again we have to describe boundedness resp. vanishing of a super function on the super upper half plane $H^{|r}$ at a cusp of $\left.\Upsilon^\# \right\backslash H$ . For 
this purpose let again $g_0 \in G$ be of the form (\ref{standardg_0}) in section 2.

\begin{lemma} \label{series} There exist series $\rund{S_n}_{n \in \nz} \in \nz^\nz$ and $\rund{D_n}_{n \in \nz} \in \rund{\rz^{r \times r}_\diag}^\nz$ such that

\begin{itemize}
\item[(i)] $\lim_{n \rightarrow \infty} S_n = + \infty$ , $\lim_{n \rightarrow \infty} D_n = 0$ , 
\item[(ii)] $\exp\rund{2 \pi i D_n} = E_0^{S_n}$ , $e^{\pi i \, \tr D_n} = \eps_0^{S_n}$ and therefore $g_0^{S_n} = \exp \chi_n$ for all $n \in \nz$ with $\chi_n := \chi_n^\diag + \chi_n^\nilp$ ,

\[
\chi_n^\diag := 2 \pi i \rund{\begin{array}{c|c}
\frac{1}{2} \, \tr D_n \, 1 & 0 \\ \hline
0 & D_n
\end{array}} \, , \, \chi_n^\nilp := \rund{\begin{array}{c|c}
\begin{array}{cc}
0 & S_n \\
0 & 0
\end{array} & 0 \\ \hline
0 & 0
\end{array}} \in \g_0 \, .
\]

\end{itemize}

\end{lemma}

{\it Proof:} simple Dirichlet argument. $\Box$ \\

Now let $\widetilde{g_0} \in_\P \G$ such that $\widetilde{g_0}^{\#'} = g_0$ .

\begin{theorem} \label{paramparab} For large $n \in \nz$ :

\item[(i)] There exist unique $\widetilde{\chi_n} \in \rund{\P \otimes \g}_0$ such that $\widetilde{\chi_n}^{\#'} = \chi_n$ and $\widetilde{g_0}^{S_n} =~\exp~\widetilde{\chi_n}$~. $s\Ad_{\widetilde{g_0}} \widetilde{\chi_n} = \widetilde{\chi_n}$ , 
and $\eckig{\widetilde{\chi_m} , \widetilde{\chi_n}} = 0$ in the Lie algebra $(\P \otimes \g)_0$ for all $m, n \in \nz$ large enough.

\item[(ii)] There exist $\P$- super automorphisms $\Omega_n$ of $H^{|r}$ such that $\Omega_n^{\#'} = \Id$ , for all $t \in \rz$

\[
\begin{array}{ccc}
\phantom{123456789} H^{|r} & \mathop{\longrightarrow}\limits^{\Omega_n} & H^{|r} \phantom{12345678} \\
\exp\rund{t \chi_n} \downarrow & \circlearrowleft & \downarrow \exp\rund{t \widetilde{\chi_n}} \\
\phantom{123456789} H^{|r} & \mathop{\longrightarrow}\limits_{\Omega_n} & H^{|r} \phantom{12345678}
\end{array} \phantom{12} \text{ , and } \phantom{12} \begin{array}{ccc}
\phantom{123} H^{|r} & \mathop{\longrightarrow}\limits^{\Omega_n} & H^{|r} \phantom{12} \\
g_0 \downarrow & \circlearrowleft & \downarrow \widetilde{g_0} \\
\phantom{123} H^{|r} & \mathop{\longrightarrow}\limits_{\Omega_n} & H^{|r} \phantom{12}
\end{array} \, .
\]

\end{theorem}

{\it Proof:} Let $n \in \nz$ be so large that $\frac{1}{2} \tr D_n$ and all the entries of $D_n$ lie in $ \, ] \, - \frac{1}{2} , \frac{1}{2} \, [ \,$~. \\

(i) For proving existence and uniqueness of $\widetilde{\chi_n}$ it suffices to show that \\
$\exp: \cz^{(2|r) \times (2|r)} \rightarrow GL(2|r, \cz)$ is a local super diffeomorphism at $\chi_n$ , and by the super inversion theorem it is even enough to show that 
$sD \exp \rund{\chi_n}$ is bijective. But since $\chi_n$ is an ordinary point of $\g_0$ and a super differential at an ordinary point involves the odd coordinates only in first order, we may without loss of generality replace the odd coordinates of 
$\cz^{(2|r) \times (2|r)}$ resp. $GL(2|r, \cz)$ by even ones and so instead show that \\
$\exp: \cz^{(2 + r) \times (2 + r)} \rightarrow GL(2 + r, \cz)$ has bijective differential at \\
$\chi_n \in \g_0 \hookrightarrow \cz^{(2 + r) \times (2 + r)}$ . We use theorem 1.7 of chapter II section 1.4 in \cite{Helga}, which says:

\begin{quote}
Let $G$ be a Lie group with Lie algebra $\g$ . The exponential mapping of the manifold $\g$ into $G$ has the differential

\[
D \exp_X = D \rund{l_{\exp X}}_e \circ \frac{1 - e^{- \ad_X}}{\ad_X}  \phantom{12345} (X \in \g) \, .
\]

As usual, $\g$ is here identified with the tangent space $\g_X$ .
\end{quote}

Hereby $e$ denotes the unit element of the Lie group $G$ , and $l_g$ denotes the left translation on $G$ with an element $g \in G$ . \\

Clearly $\ad_{\chi_n} = \ad_{\chi_n^\diag} + \ad_{\chi_n^\nilp}$ with nilpotent $\ad_{\chi_n^\nilp}$ . $\ad_{\chi_n^\diag}$ is diagonalizable and its eigenvalues are differences of the eigenvalues of $\chi_n^\diag$ and therefore 
$\in \phantom{1} ] - 2 i \pi , 2 i \pi [ \,$ . So $\frac{1 - e^{- \ad_{\chi_n}} }{\ad_{\chi_n}}$ is trigonalizable with all eigenvalues different from $0$ , which shows that $\exp$ is indeed a local super diffeomorphism at $\chi_n$ . \\

Now $s\Ad_{\widetilde{g_0}} \widetilde{\chi_n} \in \rund{\P \otimes \g}_0$ has relative body $\Ad_{g_0} \chi_n = \chi_n$ , and

\[
\exp\rund{s\Ad_{\widetilde{g_0}} \widetilde{\chi_n}} = \widetilde{g_0} \rund{\exp \widetilde{\chi_n}} \widetilde{g_0}^{- 1} = \widetilde{g_0}^{S_n} \, .
\]

Therefore by the uniqueness of $\widetilde{\chi_n}$ we see that $s\Ad_{\widetilde{g_0}} \widetilde{\chi_n} = \widetilde{\chi_n}$ , and so $\widetilde{g_0}$ commutes with all $\exp\rund{t \widetilde{\chi_n}}$ , $t \in \rz$ . Furthermore let $t \in \rz$ be arbitrary. 
Then $s\Ad_{\exp\rund{t \widetilde{\chi_m}} } \widetilde{\chi_n} \in \rund{\P \otimes \g}_0$ has relative body $\Ad_{\exp\rund{t \chi_m}} \chi_n = \chi_n$ , and

\[
\exp\rund{s\Ad_{\exp\rund{t \widetilde{\chi_m}} \widetilde{\chi_n}} } = \exp\rund{t \widetilde{\chi_m}} \rund{\exp \widetilde{\chi_n}} \exp\rund{- t \widetilde{\chi_m}} = \widetilde{g_0}^{S_n} \, .
\]

Again by the uniqueness of $\widetilde{\chi_n}$ we see that $s\Ad_{\exp\rund{t \widetilde{\chi_m}} } \widetilde{\chi_n} = \widetilde{\chi_n}$ . So $\eckig{\widetilde{\chi_m} , \widetilde{\chi_n}} = 0$ . $\Box$ \\

(ii) Take any norm $\abs{\phantom{1}}$ on the finite dimensional complex algebra $\rund{\P^\cz}^{(2|r) \times (2|r)}$ . Then there exists $C > 0$ such that $\abs{X Y} \leq C \abs{X} \abs{Y}$ for all 
$X, Y \in \rund{\P^\cz}^{(2|r) \times (2|r)}$. Clearly $\widetilde{\chi_n} \in \rund{\P \otimes \g}_0 \sqsubset \rund{\rund{\P^\cz}^{(2|r) \times (2|r)} }_0$ , and $\exp\rund{t \widetilde{\chi_n}} \in_\P \G$ , $t \in \rz$ , can be computed via ordinary exponential 
series

\[
\exp\rund{t \widetilde{\chi_n}} = \sum_{m = 0}^\infty \frac{1}{m !} t^m \widetilde{\chi_n}^m \, ,
\]

whose components are everywhere convergent power series in $t$ since \\
$\abs{\widetilde{\chi_n}^m} \leq C^{m - 1} \abs{\widetilde{\chi_n}}$ for all $m \in \nz$ . Let

\[
\exp\rund{t \widetilde{\chi_n}} = \rund{\begin{array}{cc|c}
a(t) & b(t) & \mu(t) \\
c(t) & d(t) & \nu(t) \\ \hline
\rho(t) & \sigma(t) & E(t)
\end{array}} \in \rund{\rund{\P^\cz}^{(2|r) \times (2|r)} }_0 [[t]] \, .
\]

Then by $\rund{\exp\rund{t \widetilde{\chi_n}} }^{\#'} = \exp\rund{t \chi_n}$ we see that $c(t)^{\#'} = 0$ and \\
$d(t)^{\#'} = e^{\pi i \, \tr D_n t}$ . Therefore 

\[
1 - e^{- \pi i \, \tr D_n t} \rund{c(t) i + d(t) + \nu(t) \zeta} \in \rund{\m^\cz \boxtimes \bigwedge\rund{\cz^r}}_0 [[t]]
\]

is nilpotent, more precisely its $N$-th power vanishes. Therefore all components of

\begin{eqnarray*}
&& \exp\rund{t \widetilde{\chi_n}} \rund{\begin{array}{c}
i \\ \hline
\zeta
\end{array}} = \frac{1}{c(t) i + d(t) + \nu(t) \zeta} \rund{\begin{array}{c}
a(t) z + b(t) + \mu(t) \zeta \\ \hline
\rho(t) i + \sigma(t) + E(t) \zeta
\end{array}}  \\
&& \phantom{123} = e^{- \pi i \, \tr D_n t} \sum_{m = 0}^{N - 1} \rund{1 - e^{- \pi i \, \tr D_n t} \rund{c(t) i + d(t) + \nu(t) \zeta}}^m \times \\
&& \phantom{1234567} \times \rund{\begin{array}{c}
a(t) z + b(t) + \mu(t) \zeta \\ \hline
\rho(t) i + \sigma(t) + E(t) \zeta
\end{array}} \\
&& \phantom{1} \in \eckig{\rund{\P^\cz \otimes \bigwedge\rund{\cz^r}}_0 \oplus \rund{\P^\cz \otimes \bigwedge\rund{\cz^r}}_1^{\oplus r}} [[t]]
\end{eqnarray*}

are everywhere convergent power series in $t$ . Write

\[
\exp\rund{t \widetilde{\chi_n}} \rund{\begin{array}{c}
i \\ \hline
\zeta
\end{array}} = \rund{\begin{array}{c}
f \\ \hline
\eta
\end{array}} (t, \zeta) \, .
\]

Then since $\rund{\exp\rund{t \widetilde{\chi_n}} \rund{\begin{array}{c}
i \\ \hline
\zeta
\end{array}} }^{\#'} = \exp\rund{t \chi_n} \rund{\begin{array}{c}
i \\ \hline
\zeta
\end{array}}$ we see that \\
$f (t, \zeta)^{\#'} = t + i$ and $\eta(t, \zeta)^{\#'} = \exp\rund{2 \pi i t \rund{D_n - \frac{1}{2} \tr D_n \, 1}} \zeta$ . Now define the $\P$- super automorphism $\Omega$ of $H^{|r}$ by

\[
\rund{\begin{array}{c}
f \\ \hline
\eta
\end{array}} \rund{z - i, \exp\rund{2 \pi i (z - i) \rund{\frac{1}{2} \tr D_n \, 1 - D_n}} \zeta} \, .
\]

Then $\Omega^{\#'} = \Id$ , and we prove that $\Omega$ fulfills the first commutative diagramme, in other words it transforms the action of $\exp\rund{t \chi_n}$ into the action of $\exp\rund{t \widetilde{\chi_n}}$ . Since the commutativity of the diagramme is 
equivalent to the equality of two tuples of holomorphic functions on $H$ it suffices to prove its commutativity on the non discrete subset $\rz + i \subset H$ . So let $t, u \in \rz$ . Then

\begin{eqnarray*}
&& \Omega \rund{\exp\rund{t \chi_n}\rund{\begin{array}{c}
u + i \\ \hline
\zeta
\end{array}} } = \Omega \rund{\begin{array}{c}
u + i + t \\ \hline
\exp\rund{\pi i t \rund{2 D_n - \tr D_n \, 1}} \zeta
\end{array}} \\
&& \phantom{123} = \rund{\begin{array}{c}
f \\ \hline
\eta
\end{array}} \rund{t + u, \exp\rund{\pi i u \rund{\tr D_n \, 1 - 2 D_n}} \zeta} \\
&& \phantom{123} = \exp\rund{(t + u) \widetilde{\chi_n}} \rund{\begin{array}{c}
i \\ \hline
\exp\rund{\pi i u \rund{\tr D_n \, 1 - 2 D_n}} \zeta
\end{array}} \\
&& \phantom{123} = \exp\rund{t \widetilde{\chi_n}} \exp\rund{u \widetilde{\chi_n}} \rund{\begin{array}{c}
i \\ \hline
\exp\rund{\pi i u \rund{\tr D_n \, 1 - 2 D_n}} \zeta
\end{array}} \\
&& \phantom{123} = \exp\rund{t \widetilde{\chi_n}} \rund{\begin{array}{c}
f \\ \hline
\eta
\end{array}} \rund{u, \exp\rund{\pi i u \rund{\tr D_n \, 1 - 2 D_n}} \zeta} \\
&& \phantom{123} = \exp\rund{t \widetilde{\chi_n}} \Omega \rund{\begin{array}{c}
u + i \\ \hline
\zeta
\end{array}} \, .
\end{eqnarray*}

Since finally $\widetilde{g_0}$ commutes with all $\exp\rund{t \widetilde{\chi_n}}$ , $t \in \rz$ ,

\[
\Omega_n := \frac{1}{S_n} \sum_{\sigma \in \zz \left/ S_n \zz\right.} \widetilde{g_0}^\sigma \circ \Omega \circ g_0^{- \sigma}
\]

has all the desired properties. $\Box$ \\

{\bf From now on} we will heavily use that $j\rund{g, \Z} = \Ber sD g \rund{\Z}^{\frac{1}{2 - r}}$ for all $g \in_\P \G$ , and therefore {\bf we have to assume $r \not= 2$ .}

\begin{defin} \label{parambound} \end{defin}

\begin{itemize}

\item[(i)] Let $R > 0$ and $f \in \P^\cz \boxtimes \O\rund{\{\Im z > R\}^{|r}}$ such that $f|_{\widetilde{g_0}} = f$ . Then

\[
f|_{\Omega_n} := f\rund{\Omega_n \rund{\Z}} \rund{\Ber sD \Omega_n}^{\frac{k}{2 - r}}
\]

is invariant under $|_{g_0}$ . $f$ is called bounded (vanishing) at $i \infty$ iff $f|_{\Omega_n}$ is bounded (vanishing) at $i \infty$ in the sense of definition \ref{bound} (i) for almost all $n \in \nz$ .

\item[(ii)] Let $z_0 \in \partial_{\pz^1} H$ and $\gamma \in_\P \G$ such that $\gamma^\# \in N^{z_0} \setminus \{\id\}$ . Let $U \subset H$ be an open and $\gamma^\#$-invariant neighbourhood of $z_0$ and $f \in \P^\cz \otimes \O\rund{U^{|r}}$ such that 
$f|_\gamma = f$ . Take some $g \in_\P \G$ such that $g^\# i \infty = z_0$ and either \\
$g_0 := \rund{g^{\#'}}^{- 1} \gamma^{\#'} g^{\#'}$ or $g_0 := \rund{g^{\#'}}^{- 1} \rund{\gamma^{\#'}}^{- 1} g^{\#'}$ is of the form (\ref{standardg_0}) (in fact we always find ordinary elements in $G$ providing this). Then \\
$\widetilde{g_0} := g^{- 1} \gamma g \in_\P \G$ resp. $\widetilde{g_0} := g^{- 1} \gamma^{- 1} g \in_\P \G$ has relative body $g_0$ , and $f|_g$ is invariant under $|_{\widetilde{g_0}}$ . Now $f$ is called 
bounded (vanishing) at $z_0$ iff $f|_g$ is bounded (vanishing) at $i \infty$ .

\end{itemize}

Observe that all powers $\rund{\Ber sD \Omega_n}^u$ , $u \in \rz$ , are well defined since $\Omega^{\#'} = \Id$ and so $\rund{\Ber sD \Omega_n}^{\#'} = 1$ . \\

Of course we have to prove well-definedness in definition \ref{parambound}, which is not at all trivial. For $D = \rund{\begin{array}{ccc}
\delta_1 & & 0 \\
 & \ddots & \\
0 & & \delta_r
\end{array}} \in \rz^{r \times r}_\diag$ and $I \in \wp(r)$ let $\tr_I D := \sum_{i \in I} \delta_i$ . Then $\det_I \exp\rund{2 \pi i D} = e^{2 \pi i \tr_I D}$ . Let us start with the independence of (i) of the choices of the $\P$- super automorphisms $\Omega_n$ of 
$H^{|r}$~:

\begin{quote}
Let $I \in \wp(r)$ . If $\eps_0^{- k - \abs{I}} \det_I E_0 \not= 1$ then

\[
\Delta_I := \min\schweif{\abs{\mu} \, \left| \, \mu \in \rz \, , \, e^{2 \pi i \mu} = \eps_0^{k + \abs{I}} \det\nolimits_I E_0^{- 1} \right.} > 0 \, .
\]

Clearly $\tr_I D_n - \frac{\abs{I} + k}{2} \, \tr D_n \rightarrow 0$ for $n \rightarrow \infty$ . The independence is shown by the following lemma:
\end{quote}

\begin{lemma} Assume that $n \in \nz$ is so large that for all $I \in \wp(r)$

\[
\abs{\tr_I D_n - \frac{\abs{I} + k}{2} \, \tr D_n} < 
\left\{\begin{array}{ll}
1 & \text{ if } \eps_0^{- k - \abs{I}} \det_I E_0 = 1 \\
S_n \Delta_I & \text{ if } \eps_0^{- k - \abs{I}} \det_I E_0 \not= 1
\end{array}\right. \, ,
\]

and let $\Omega$ be a $\P$- super automorphism of $H^{|r}$ having $\Omega^{\#'} = \Id$ and commuting with all $\exp\rund{t \chi_n}$ , $t \in \rz$ . Let $f \in \O\rund{\{\Im z > R\}^{|r}}$ be invariant under $|_{g_0}$ . Then if $f$ is bounded (vanishing) at 
$i \infty$ so is $f|_\Omega$ .
\end{lemma}

{\it Proof:} Let $\Xi$ be the super automorphism of $H^{|r}$ given by $\rund{\begin{array}{c}
z \\ \hline
\exp\rund{\pi i \frac{z}{S_n} \rund{2 D_n - \tr D_n \, 1}} \zeta
\end{array}}$ . Then $\Xi^\# = \id$ , and straight forward computations show that $\Ber sD \Xi = e^{\frac{\pi i (r - 2) \, \tr D_n}{S_n} z}$ ,

\[
\begin{array}{ccc}
\phantom{123456789} H^{|r} & \mathop{\longrightarrow}\limits^{\Xi} & H^{|r} \phantom{12345678} \\
\exp\rund{t \chi_n^\nilp} \downarrow & \circlearrowleft & \downarrow \exp\rund{t \chi_n} \\
\phantom{123456789} H^{|r} & \mathop{\longrightarrow}\limits_{\Xi} & H^{|r} \phantom{12345678}
\end{array} \, ,
\]

and $f|_\Xi = f|_\Xi\rund{\begin{array}{c}
z + 1 \\ \hline
\zeta
\end{array}} = \left.f|_\Xi \phantom{\frac{1}{1}} \right|_{\exp \chi_n^\nilp}$ , where

\[
f|_\Xi := f\rund{\Xi \rund{\Z}} e^{\frac{\pi i k \, \tr D_n}{S_n} z} \, .
\]

First we show that $f$ is bounded (vanishing) at $i \infty$ iff $f|_\Xi$ is bounded (vanishing) at $i \infty$ .

\begin{quote}
Since $|_\Xi$ respects the splitting $f = \sum_{I \in \wp(r)} f_I \zeta^I$ we may assume without restriction that $f = f_I \zeta^I$ for some $I \in \wp(r)$ and \\
$f_I \in \O\rund{\{\Im z > R\}}$ . Then

\[
f|_{\Xi} = e^{\pi i \frac{2 \tr_I D_n - (\abs{I} + k) \, \tr D_n}{S_n} \, z} f \, .
\]

{\it First case:} $\eps_0^{- k - \abs{I}} \det_I E_0 = 1$ . Then $\tr_I D_n - \frac{\abs{I} + k}{2} \tr D_n < 1$ by the assumption on $n$ , and so since on the other hand 
$e^{\pi i \rund{2 \tr_I D_n - (\abs{I} + k) \tr D_n}} = \rund{\eps_0^{- k - \abs{I}} \det_I E_0}^{S_n} = 1$ we see that $\tr_I D_n - \frac{\abs{I} + k}{2} \, \tr D_n = 0$ , and so $f|_{\Xi} = f$ . \\

{\it Second case:} $\eps_0^{- k - \abs{I}} \det_I E_0 \not= 1$ . $f|_{g_0} = f$ implies \\
$f_I = \eps_0^{- k - \abs{I}} \det_I E_0 \, f_I(z + 1)$ , and so by Fourier decomposition we may assume without loss of generality that $f_I = e^{2 \pi i \mu z}$ for some $\mu \in \rz$ , $e^{2 \pi i \mu} = \eps_0^{k + \abs{I}} \det_I E_0^{- 1}$ . So

\[
f|_{\Xi} = e^{\pi i \rund{2 \mu + \frac{2 \tr_I D_n - (\abs{I} + k) \tr D_n}{S_n}} z} \zeta^I \, .
\]

Assume $f$ bounded at $i \infty$ . Then $\mu \geq 0$ , and so $\mu \geq \Delta_I$~. Since $\tr_I D_n - \frac{\abs{I} + k}{2} \, \tr D_n < S_n \Delta_I$ by assumption on $n$ , we have $2 \mu + \frac{2 \tr_I D_n - (\abs{I} + k) \, \tr D_n}{S_n} > 0$ and so $f|_\Xi$ is 
in fact even vanishing at $i \infty$ .

Conversely assume $f$ {\bf not} vanishing at $i \infty$ . Then $\mu \leq 0$ , so $\mu \leq - \Delta_I$ and therefore $2 \mu + \frac{2 \tr_I D_n - (\abs{I} + k) \, \tr D_n}{S_n} < 0$ . We see that in this case $f|_\Xi$ is even {\bf not} bounded at $i \infty$ .
\end{quote}

So replacing $f$ by $f|_\Xi$ and $\Omega$ by $\, \Xi^{- 1} \circ \Omega \circ \Xi \,$ we may assume without loss of generality that $f|_{\exp \chi_n^\nilp} = f$ , $\Omega^{\#'} = \Id$ and $\Omega$ commutes with all $\exp\rund{t \chi_n^\nilp}$ . A simple 
computation shows that then $\Omega$ must be of the form

\[
\rund{\begin{array}{c}
z + \sum_{J \in \wp(r)} a_J \zeta^J \\ \hline
\zeta + \sum_{J \in \wp(r)} b_J \zeta^J
\end{array}} \, ,
\]

all $a_J \in \m^\cz$ , $b_J \in \rund{\m^\cz}^{\oplus r}$ of suitable parity, and therefore \\
$\Ber sD \Omega = 1 + \sum_{J \in \wp(r)} c_J \zeta^J$ with some $c_J \in \m^\cz$ . So if we assume without loss of generality that $f = f_I \zeta^I$ , $I \in \wp(r)$ , $f_I \in \O\rund{\{\Im z > R\}}$ , we obtain

\begin{eqnarray*}
f|_\Omega &=& \sum_{m = 0}^{N - 1} f_I^{(m)}(z) \rund{\sum_{J \in \wp(r)} a_J \zeta^J}^m \rund{\zeta + \sum_{J \in \wp(r)} b_J \zeta^J}^I \times \\
&& \phantom{123} \times \rund{1 + \sum_{J \in \wp(r)} c_J \zeta^J}^{\frac{k}{2 - r}} \, ,
\end{eqnarray*}

which is a linear combination over $\P^\cz$ of expressions $f_I^{(m)}(z) \zeta^K$ , $m \in \{0, \dots, N - 1\}$ , $K \in \wp(r)$ . Therefore if $f$ is bounded (vanishing) at $i \infty$ then so is $f|_\Omega$ . $\Box$ \\

The following lemma is of independent interest but will in particular show that (i) is independent of the choice of the series $\rund{S_n}_{n \in \nz}$ and $\rund{D_n}_{n \in \nz}$ :

\begin{lemma} \label{onlyone} There exists $n_0 \in \nz$ such that for all $n \geq n_0$ : $f$ is bounded (vanishing) at $i \infty$ iff $f|_{\Omega_n}$ is bounded (vanishing) at $i \infty$ in the sense of definition \ref{bound} (i).
\end{lemma}

{\it Proof:} We just have to show that for large $m, n \in \nz$ we can find a common $\Omega_m = \Omega_n$ . For this purpose let $\Omega_n$ be given by theorem \ref{paramparab} (ii). Then since all $\exp\rund{t \widetilde{\chi_m}}$ and 
$\exp\rund{u \widetilde{\chi_n}}$ , $t, u \in \rz$ commute, we see that

\[
\Omega_m' := \Omega_n' := \int_{\rz \left/ \zz\right.} \exp\rund{2 \pi i \sigma \widetilde{\chi_m}} \circ \Omega_n \circ \exp\rund{- 2 \pi i \sigma \chi_m} d \sigma
\]

fulfills at the same time all the desired properties of both $\Omega_m$ and $\Omega_n$ in theorem \ref{paramparab} (ii). $\Box$ \\

Now let us show the independence of (ii) of the choice of $g \in_\P \G$ :

\begin{quote}
Let $g \in_\P G$ such that $g_0' := \rund{g^{\#'}}^{- 1} g_0 \, g^{\#'}$ is again of the form (\ref{standardg_0}) with some $\eps_0' \in U(1)$ , $E_0' \in U(r)$ diagonal, $\eps_0'^2 = \det E_0'$ . Then $\eps_0' = \eps_0$ and $E_0' = P E_0 P^{- 1}$ with 
some permutation matrix $P \in U(r)$ . $f|_g$ is invariant under $|_{\widetilde{g_0}'}$ , $\widetilde{g_0}' := g^{- 1} \widetilde{g_0} g$ , and we have to prove that if $f$ is bounded (vanishing) at $i \infty$ then so is $f|_g$ . Let the series 
$\rund{S_n}_{n \in \nz}$ and $\rund{D_n}_{n \in \nz}$ be given by lemma \ref{series} with respect to $g_0$ . Then the series $\rund{S_n}_{n \in \nz}$ and $\rund{D_n'}_{n \in \nz}$~, $D_n' := P D_n P^{- 1}$ , and the resulting $\chi_n' \in \g_0$ , $n \in \nz$ , 
fulfill all the desired properties of lemma \ref{series} with respect to $g_0'$ instead of~$g_0$~.

\begin{lemma} Let $n \in \nz$ be so large that all the entries of $D_n$ lie in $\, \left] - \frac{1}{2} , \frac{1}{2} \right[ \,$ . Then $\Ad_{\rund{g^{\#'}}^{- 1}} \chi_n = \chi_n'$ .
\end{lemma}

{\it Proof:} Let $E \in U(r)$ be the lower right corner of $g^{\#'}$ . Then obviously $E P$ commutes with $E_0$ . So since $\exp\rund{2 \pi i D_n}$ is the lower right corner of $g_0^{S_n}$ we see that $E P$ stabilizes all eigenspaces of 
$\exp\rund{2 \pi i D_n}$ . But all the eigenvalues of $D_n$ lie in $\, \left] - \frac{1}{2} , \frac{1}{2} \right[ \,$ . Therefore the eigenspaces of $D_n$ are the same as the ones of $\exp\rund{2 \pi i D_n}$ . So $E P$ even commutes with $D_n$ . This implies 
$\Ad_{g^{\#'} h} \chi_n = \chi_n$ , and so $\Ad_{\rund{g^{\#'}}^{- 1}} \chi_n = \Ad_h \chi_n = \chi_n'$~.~$\Box$ \\

Therefore for large $n \in \nz$ : $\widetilde{\chi_n}' := s\Ad_{\rund{g^{\#'}}^{- 1}} \widetilde{\chi_n}$ are the unique elements of  $(\P \otimes \g)_0$ given by theorem \ref{paramparab} (i) with respect to $\widetilde{g_0}'$ and $\chi_n'$ instead 
of $\widetilde{g_0}$ resp. $\chi_n$ , and $\Omega_n' := g^{- 1} \circ \Omega_n \circ g^{\#'}$ fulfill all the desired properties in theorem \ref{paramparab} (ii) with respect to $\widetilde{g_0}'$ and $\chi_n'$ instead of $\widetilde{g_0}$ resp. $\chi_n$ . So we 
have to show that $f|_{\Omega_n}$ bounded (vanishing) at $i \infty$ implies $\left.f|_g \right|_{\Omega_n'} = \left.f|_{\Omega_n} \right|_{g^{\#'}}$ bounded (vanishing) at $i \infty$ , which has already been proven for the well-definedness of definition 
\ref{bound}.
\end{quote}

Of course (ii) still depends on the choice of $\gamma$ . However, let us show that (ii) is invariant under replacing $\gamma \in_\P \G$ be some power $\gamma^m$ , $m \in \zz \setminus \{0\}$ :

\begin{quote}
Without loss of generality we may assume that $m \in \nz \setminus \{0\}$ and $\widetilde{g_0} = \gamma^m$ . Let $g := \rund{\begin{array}{c|c}
\begin{array}{cc}
\frac{1}{\sqrt{m}} & 0 \\
0 & \sqrt{m}
\end{array} & 0 \\ \hline
0 & 1
\end{array}} \in G$ . Then \\
$g_0' := g^{- 1} \gamma^{\#'} g$ is again of the form \ref{standardg_0} with some $\eps_0' \in U(1)$~, \\
$E_0' \in U(r)$ diagonal, $\eps_0'^2 = \det E_0'$ such that $\eps_0'^m = \eps_0$ and \\
$E_0'^m = E_0$ . Let the series $\rund{S_n}_{n \in \nz}$ and $\rund{D_n}_{n \in \nz}$ be given by lemma \ref{series} with respect to $g_0$ . Then the series $\rund{S_n'}_{n \in \nz}$ given by $S_n' := m S_n$ and $\rund{D_n}_{n \in \nz}$ and the resulting 
$\chi_n'$ , $n \in \nz$ , fulfill lemma \ref{series} with respect to $g_0'$ .

Furthermore let $\widetilde{\chi_n}$ and $\Omega_n$ be given by theorem \ref{paramparab} with respect to $\widetilde{g_0}$ and $\widetilde{\chi_n}'$ with respect to $\widetilde{g_0}'$ . Then we obtain 
$\widetilde{\chi_n}' = s\Ad_{g^{- 1}} \widetilde{\chi_n}$~, and $\Omega_n' := g^{- 1} \circ \Omega_n \circ g$ fulfills all the desired properties in theorem \ref{paramparab} (ii) with respect to $\widetilde{g_0}'$ . So we have to show that 
$\left.f|_g \right|_{\Omega_n'} = \left.f|_{\Omega_n} \right|_g$ is bounded (vanishing) at $i \infty$ iff so is $f|_{\Omega_n}$ , which is quite obvious.
\end{quote}

Let $\Upsilon$ be a $\P$-lattice of $\G$ and $\Upsilon_0$ denote the kernel of the body map $\Upsilon \rightarrow \Aut H$ or equivalently the preimage of $\rund{\Upsilon^{\#'}}_0$ in $\Upsilon$ under ${}^{\#'}$ . \\

Assume $\gamma \in \Upsilon$ . Then definition \ref{parambound} (ii) is even invariant under replacing $\gamma$ by another element $\eta \in \Upsilon$ having $\eta^\# = \gamma^\#$ in the case where $f$ is also invariant under $|_\eta$ , which is a trivial 
consequence of the invariance of (ii) under replacing $\gamma \in_\P \G$ be some power $\gamma^m$ and the following lemma.

\begin{lemma}
Let $\gamma, \eta \in \Upsilon$ having $\gamma^\# = \eta^\#$ . Then there exists some \\
$m \in \nz \setminus \{0\}$ such that $\gamma^m = \eta^m$ .
\end{lemma}

{\it Proof:} Clearly $\rund{\gamma^l \eta^{- l}}^\# = \id$ and so $\gamma^l \eta^{- l} \in \Upsilon_0$ for all $l \in \nz$~. But $\Upsilon_0$ is finite, so there exist $l, l' \in \nz$ such that $l > l'$ and $\gamma^{l} \eta^{- l} = \gamma^{l'} \eta^{- l'}$ . Taking $m := l - l'$ 
yields $\gamma^m = \eta^m$ . $\Box$ \\

Now we are ready for giving the definition of super automorphic and super cusp forms for the $\P$-lattice $\Upsilon$ .

\begin{defin}[super automorphic and super cusp forms for $\Upsilon$ ] \end{defin}

Let $f \in \P^\cz \otimes \O\rund{H^{|r}}$ . $f$ is called a super automorphic (cusp) form for $\Upsilon$ of weight $k$ iff

\begin{itemize}
\item[(i)] $f|_\gamma = f$ for all $\gamma \in \Upsilon$ ,
\item[(ii)] $f$ is bounded (vanishing) at all cusps $\overline{z_0} \in {\Upsilon^{\#}} \left\backslash \partial_{\pz^1} H^{\phantom{1}} \right.$ of $\left.\Upsilon^\# \right\backslash H$ in the sense of definition \ref{parambound}.
\end{itemize}

The $\zz_2$-graded $\P^\cz$-module of super automorphic (cusp) forms for $\Upsilon$ of weight $k$ is denoted by $sM_k(\Upsilon)$ (resp. $sS_k(\Upsilon)$ ). In general these spaces do {\bf not} have a canonical $\zz$-grading! \\

As a trivial observation let us remark that if $g \in_\P \G$ then \\
$sM_k(\Upsilon) \mathop{\rightarrow}\limits^\sim sM_k\rund{g^{- 1} \Upsilon g} \, , \, f \mapsto f|_g$ is a graded isomorphism mapping $sS_k(\Upsilon)$ to $sS_k\rund{g^{- 1} \Upsilon g}$ . Furthermore $\rund{f|_g}^{\#'} = \left.f^{\#'}\right|_{g^{\#'}}$ for all \\
$f \in \P^\cz \boxtimes \O\rund{g^\# U^{|r}}$ , $U \subset H$ open, and $g \in_\P \G$ , so in particular ${}^{\#'}$ restricts to a linear map

\[
{}^{\#'} : sM_k(\Upsilon) \rightarrow sM_k\rund{\Upsilon^{\#'}}
\]

mapping $sS_k(\Upsilon)$ to $sS_k\rund{\Upsilon^{\#'}}$ . \\

From now on let $\Gamma \sqsubset G$ be a lattice and $k_2 \in \zz$ be given by lemma \ref{H1vanish}. We may assume $k_2$ to be independent of the choice of $\rho \in \{0, \dots, r\}$ by taking the maximum over all $\rho$ .

\begin{theorem}[main theorem] \label{parammain} For any $\P$-lattice $\Upsilon$ of $\G$ with relative body $\Gamma$ and weight $k \geq k_2$ we have $\zz_2$-graded $\P^\cz$-module isomorphisms

\[
\begin{array}{ccc}
sS_k(\Upsilon) & \simeq & \P^\cz \otimes sS_k(\Gamma) \\
\cap & \circlearrowleft & \cap \\
sM_k(\Upsilon) & \simeq & \P^\cz \otimes sM_k(\Gamma) \\
{}_{{}^{\#'}} \searrow & \circlearrowleft & \swarrow_{{}^{\#'} \otimes \id} \\
 & sM_k(\Gamma) & \phantom{12345678901234567890} .
\end{array}
\]

\end{theorem}

We will show that this is a special case of the situation discussed in section 4 with $K := \cz$ and $\P^\cz$ instead of $\P$ . \\

\begin{quote}
Let us again briefly discuss example \ref{exautfo} (ii): In both cases $\sdim H^1\rund{\Gamma, \g} = (1, 2)$ . Let $\Upsilon$ be a $\P$-lattice of $G$ with $\Upsilon^{\#'} = \Gamma$ . \\

In both cases we have indeed a $\zz_2$-graded $\P^\cz$-module isomorphism 

\[
\begin{array}{ccc}
sM_k(\Upsilon) & \simeq & \P^\cz \otimes sM_k(\Gamma) \\
{}_{{}^{\#'}} \searrow & \circlearrowleft & \swarrow_{{}^{\#'} \otimes \id} \\
 & sM_k(\Gamma) & \phantom{12345678901234567890} ,
\end{array}
\]

$k = 1$ in the first and $k = 0$ in the second case. This is evident applying the proof of theorem \ref{parammain}, in particular lemma \ref{local} and lemma \ref{adapt}, to this special 
situation using the fact that \\
$H^1\rund{E_1^0} = H^1\rund{E_1^1} = 0$ in the first and $H^1\rund{E_0^0} = H^1\rund{E_0^1} = 0$ in the second case. In particular we see that there exists in the 
first case $\widetilde{\eta^2} \in sM_1(\Upsilon)$ even with $\widetilde{\eta^2}^{\#'} = \eta^2$ and in the second case $\widetilde{\eta^2 \zeta} \in sM_0(\Upsilon)$ odd with 
$\widetilde{\eta^2 \zeta}^{\#'} = \eta^2\zeta$ .
\end{quote}

Now define the sheaves $\F_k \hookrightarrow \E_k$ of $\zz_2$-graded $\P^\cz$-modules on $X$ as

\begin{eqnarray*}
&& \E_k(U) := \left\{f \in \P^\cz \otimes \O\rund{\pi_X^{- 1}(U)^{|r}} \, \right| \, f|_\gamma = f \text{ for all } \gamma \in \Upsilon \, , \\
&& \left. \phantom{1234 \pi_X^{- 1} \P^\cz} f \text{ bounded at all cusps } \overline{z_0} \in U \text{ of } \Gamma^\# \backslash H\right\}
\end{eqnarray*}

and

\begin{eqnarray*}
&& \F_k(U) := \left\{f \in \P^\cz \otimes \O\rund{\pi_X^{- 1}(U)^{|r}} \, \right| \, f|_\gamma = f \text{ for all } \gamma \in \Upsilon \, , \\
&& \left. \phantom{1234 \pi_X^{- 1} \P^\cz} f \text{ vanishing at all cusps } \overline{z_0} \in U \text{ of } \Gamma^\# \backslash H\right\}
\end{eqnarray*}

for all $U \subset X$ open. Recall that $\pi_X : H \rightarrow \left.\Gamma^\# \right\backslash H \hookrightarrow X$ denotes the canonical projection. Clearly $sM_k(\Upsilon) = \E_k(X)$ and $sS_k(\Upsilon) = \F_k(X)$ . ${}^{\#'}$ induces a graded sheaf 
projection ${}^{\#'} : \E_k \rightarrow \left.\E_k \right/ \m \E_k \simeq \Gamma^\hol\rund{\diamondsuit, E_k}$ restricting to a sheaf projection $\F_k \rightarrow \left.\F_k \right/ \m \F_k \simeq \Gamma^\hol(\diamondsuit, F_k)$ , where \\
$E_k := \bigoplus_{\rho = 0}^r E_k^\rho$ and $F_k := \bigoplus_{\rho = 0}^r F_k^\rho \rightarrow X$ denote the holomorphic vector bundles from section 2.

\begin{lemma} \label{local} Locally we have $\zz_2$-graded $\P^\cz$-module isomorphisms

\[
\begin{array}{ccc}
\F_k & \simeq & \P^\cz \otimes \Gamma^\hol\rund{\diamondsuit, F_k} \\
\cap & \circlearrowleft & \cap \\
\E_k & \simeq & \P^\cz \otimes \Gamma^\hol\rund{\diamondsuit, E_k} \\
{}_{{}^{\#'}} \searrow & \circlearrowleft & \swarrow_{{}^{\#'} \otimes \id} \\
 & \Gamma^\hol\rund{\diamondsuit, E_k} & \phantom{12345678901234567890} .
\end{array}
\]

\end{lemma}

{\it Proof:} {\it First case:} $z_0 \in H$ . Then $\Gamma^{z_0} := \schweif{\gamma \in \Gamma \, \left| \, \gamma^\# z_0 = z_0\right.} \sqsubset \Gamma$ is a finite subgroup. Since the action of $\Aut H$ is proper, there exists an open and 
$\rund{\Gamma^{z_0}}^\#$-invariant neighbourhood $U \subset H$ of $z_0$ such that $\pi_X(U)$ is an open neighbourhood of $\overline{z_0}$ in $X$ and $\pi_X$ induces a biholomorphic map 
$\left.\Gamma^{z_0} \right\backslash U \simeq \pi_X(U)$ . 
Via this biholomorphic map we obtain a graded sheaf homomorphism $\varphi$ from $\left.\E_k\right|_{\pi_X(U)} = \left.\F_k\right|_{\pi_X(U)}$ to \\
$\P^\cz \otimes \Gamma^\hol\rund{\diamondsuit, \left.E_k\right|_{\pi_X(U)} } = \P^\cz \otimes \Gamma^\hol\rund{\diamondsuit, \left.F_k\right|_{\pi_X(U)} }$ respecting ${}^{\#'}$ given by

\[
\varphi_V : \E_k(V) \rightarrow \P^\cz \otimes \Gamma^\hol\rund{V, E_k} \, , \, f \mapsto \frac{1}{\abs{\Gamma^{z_0}} } \sum_{\gamma \in \Gamma^{z_0}} f|_\gamma
\]

for all $V \subset \pi_X(U)$ open, where we use the canonical embeddings $\E_k(V) , \P^\cz \otimes \Gamma^\hol\rund{V, E_k} \hookrightarrow \P^\cz \boxtimes \O\rund{\rund{\pi_X^{- 1}(V) \cap U}^{|r}}$ . \\

For proving injectivity of $\varphi_V$ , $V \subset U$ open, let $f \in \E_k(V) \setminus \{0\}$ and $n \in \nz$ be maximal such that $f \in \rund{\m^n}^\cz \otimes \O\rund{V^{|r}}$ . Then in particular 
$\varphi_V(f) \equiv f \not\equiv 0 \mod \rund{\m^{n + 1}}^\cz \otimes \O\rund{V^{|r}}$ , which shows that $\varphi_V(f) \not= 0$ . \\

Now assume that $\varphi_V$ is not surjective for some $V \subset \pi_X(U)$ open and let $n \in \nz$ be maximal with the property that there exists \\
$h \in \rund{\m^n}^\cz \otimes \Gamma^\hol\rund{V, E_k} \setminus \Im \varphi_V$ . Then with \\
$\Upsilon^{z_0} := \schweif{\gamma \in \Upsilon \, \left| \, \gamma^{\#'} \in \Gamma^{z_0}\right.} \sqsubset \Upsilon$ we see that

\[
f := \frac{1}{\abs{\Gamma^{z_0}} } \sum_{\gamma \in \Upsilon^{z_0}} h|_\gamma \in \rund{\m^n}^\cz \E_k\rund{\pi_X(V)}
\]

and $f \equiv h \mod \rund{\m^{n+ 1}}^\cz \otimes \O\rund{V^{|r}}$ . Therefore $h - \varphi_V(f) \equiv 0 \mod \rund{\m^{n + 1}}^\cz \otimes \Gamma^\hol\rund{V, E_k}$ . So maximality of $n$ implies \\
$h - \varphi_V(f) \in \Im \varphi_V$~, which is a contradiction to the linearity of $\varphi_V$~. \\

{\it Second case:} $z_0 \in \partial_{\pz^1} H$ cusp of $\left.\Gamma^\# \right\backslash H$ . Let $\gamma \in \Upsilon$ such that $\gamma^\#$ generates $N^{z_0} \cap \Gamma^\#$ . Again there exists an open and $\gamma^\#$-invariant neighbourhood 
$U \subset H$ of $z_0 = i \infty$ such that $\pi_X(U) \cup \schweif{\overline{i \infty}}$ is an open neighbourhood of $\overline{i \infty}$ in $X$ , and $\pi_X$ induces a biholomorphic map $\left.\spitz{\gamma^\#} \right\backslash U \simeq \pi_X(U)$ . \\

After applying $|_g$ to the sections of both the sheaf $\E_k$ and the vector bundle $E_k$ with a suitable $g \in G$ we may assume without loss of generality that $z_0 = i \infty$ and $g_0 := \gamma^{\#'} \in \Gamma$ is of the form (\ref{standardg_0}). 
Define $\widetilde{g_0} := \gamma$ .
Let the $\P$-isomorphisms $\Omega_n$ , $n \in \nz$ large, be given by theorem \ref{paramparab} (ii) and $n_0 \in \nz$ be given by lemma \ref{onlyone}. Again we obtain a graded sheaf homomorphism $\psi$ from 
$\left.\E_k\right|_{\pi_X(U) \cup \schweif{\overline{i \infty}} }$ to $\P^\cz \otimes \Gamma^\hol\rund{\diamondsuit, \left.E_k\right|_{\pi_X(U) \cup \schweif{\overline{i \infty}} }}$ respecting ${}^{\#'}$ and mapping 
$\left.\F_k\right|_{\pi_X(U) \cup \schweif{\overline{i \infty}} }$ to $\P^\cz \otimes \Gamma^\hol\rund{\diamondsuit, \left.F_k\right|_{\pi_X(U)  \cup \schweif{\overline{i \infty}} }}$ given by

\[
\psi_V : \E_k(V) \rightarrow \P^\cz \otimes \Gamma^\hol\rund{V, E_k} \, , \, f \mapsto \frac{1}{\abs{\Gamma_0}} \sum_{\eta \in \Gamma_0} \left.f|_{\Omega_{n_0}} \right|_\eta
\]

for all $V \subset \pi_X(U) \cup \schweif{\overline{i \infty}}$ open, where again we use the canonical embeddings $\E_k(V) , \P^\cz \otimes \Gamma^\hol\rund{V, E_k} \hookrightarrow \P^\cz \boxtimes \O\rund{\rund{\pi_X^{- 1}(V) \cap U}^{|r}}$ . Injectivity of all 
$\psi_V$ is proven in a similar way to the case $z_0 \in H$ . \\

Again assume that $\psi_V$ is not surjective (resp. the preimage of \\
$\P^\cz \otimes \Gamma^\hol\rund{V, F_k}$ under $\psi_V$ does not lie in $\F_k(V)$ ) for some \\
$V \subset \pi_X(U) \cup \schweif{\overline{i \infty}}$ open, and let $n \in \nz$ be maximal with the property that there exists $h \in \rund{\m^n}^\cz \otimes \Gamma^\hol\rund{V, E_k} \setminus \Im \psi_V$ \\
(resp. $h \in \rund{\m^n}^\cz \otimes \Gamma^\hol\rund{V, F_k} \setminus \psi_V \F_k(V)$ ). Then 

\[
f := \frac{1}{\abs{\Gamma_0}} \sum_{\eta \in \Upsilon_0} \left.h|_{\Omega_{n_0}^{- 1}} \right|_\eta \in \rund{\m^n}^\cz \E_k(V) \, ,
\]

where $h|_{\Omega_{n_0}^{- 1}} := h\rund{\Omega_n^{- 1} \Z} \rund{\Ber sD \, \Omega_n^{- 1}}^{\frac{k}{2 - r}}$ \, .

\begin{quote}
Let us show that indeed $f$ is bounded (vanishing) at $i \infty$ if $\overline{i \infty} \in~V$ and $h \in \rund{\m^n}^\cz \otimes \Gamma^\hol\rund{V, E_k}$ \\
(resp. $h \in \rund{\m^n}^\cz \otimes \Gamma^\hol\rund{V, F_k}$ ): By lemma \ref{onlyone} clearly $h|_{\Omega_{n_0}^{- 1}}$ , which is invariant under $|_\gamma$ , is bounded (vanishing) at $i \infty$ . Now let $\eta \in \Upsilon_0$ . Then $\eta^\# = \id$ and 
so $\eta^{- 1} \gamma \eta \in \Upsilon$ has body $\gamma^\#$ . So $\left.h|_{\Omega_{n_0}^{- 1}} \right|_\eta$ , which is invariant under $|_{\eta^{- 1} \gamma \eta}$ , is clearly bounded (vanishing) at $i \infty$ .
\end{quote}

Again $f \equiv h \mod \rund{\m^{n+ 1}}^\cz \otimes \O\rund{V^{|r}}$ . Therefore $h - \psi_V(f) \equiv 0 \mod \rund{\m^{n + 1}}^\cz \otimes \Gamma^\hol\rund{V, E_k}$ . So again $h - \psi_V(f) \in \Im \psi_V$ (resp. $h - \psi_V(f) \in \psi_V \F_k(V)$ ), which is a 
contradiction to the linearity of $\psi_V$ . Therefore all $\psi_V$ , $V \subset \pi_X(U)$ open, are surjective. $\Box$ \\

{\it Finally for proving theorem \ref{parammain}:} $H^1\rund{E_k} , H^1\rund{F_k} = 0$ by the choice of $k$ . Let $\rund{f_1, \dots, f_d}$ be a graded basis of $sM_k(\Gamma) =~H^0\rund{E_k}$ such that $\rund{f_1, \dots, f_{d'}}$ is a basis of 
$sS_k(\Gamma) = H^0\rund{F_k}$ . Then by lemma \ref{adapt} there exist $\widetilde{f_1}, \dots, \widetilde{f_{d'}} \in \F_k(X) = sS_k(\Upsilon)$ , $\widetilde{f_{d' + 1}}, \dots, \widetilde{f_d} \in sM_k(\Upsilon)$ such that ${\widetilde{f_\delta}}^{\#'} = f_\delta$ 
for all $\delta = 1, \dots, d$ . Since ${}^{\#'}$ is graded, after applying the projections onto the even resp. odd summand of $sM_k(\Upsilon)$ to $\widetilde{f_\delta}$ , $\delta = 1, \dots, d$ , we may assume without restriction that $\widetilde{f_\delta}$ is 
graded of the same parity as $f_\delta$ , $\delta = 1, \dots, d$~. Therefore the $\P^\cz$-module isomorphism 

\[
sM_k(\Upsilon) \rightarrow \P^\cz \otimes sM_k(\Gamma)
\]

given by the assignment $\widetilde{f_\delta} \mapsto f_\delta$ is graded and has all the desired properties.~$\Box$

\end{document}